\documentclass{jT}
\usepackage{amsmath,amsfonts,amsthm,amssymb,amscd,url,stmaryrd}
\def\classification#1{\def\@class{#1}}
\classification{\null}

\newcommand{\st}{s.t.}

\DeclareFontFamily{OT1}{rsfs}{}
\DeclareFontShape{OT1}{rsfs}{n}{it}{<-> rsfs10}{}
\DeclareMathAlphabet{\mathscr}{OT1}{rsfs}{n}{it}
\DeclareMathOperator{\sgn}{sgn}

\DeclareMathOperator{\mo}{\,mod}

\DeclareMathOperator{\Prob}{\mathbb{P}}

\DeclareMathOperator{\Disc}{Disc}

\DeclareMathOperator{\li}{li}
\DeclareMathOperator{\rnk}{rank}
\DeclareMathOperator{\Div}{Div}

\DeclareMathOperator{\Frob}{Frob}
\DeclareMathOperator{\Gal}{Gal}

\DeclareMathOperator{\Cl}{Cl}

\newtheorem{theorem}{Theorem}[section]

\newtheorem{thm}[theorem]{Theorem}
\newtheorem{prop}[theorem]{Proposition}
\newtheorem{cor}[theorem]{Corollary}
\newtheorem{lem}[theorem]{Lemma}
\newtheorem*{main}{Main Theorem}

\theoremstyle{definition}
\newtheorem{defn}{Definition}[section]
\newtheorem*{remark}{Remark}
\newtheorem*{example}{Example}

\numberwithin{equation}{section}
\begin{document}
\title[Power-free values, large deviations, and integer points
on irrational curves]{Power-free values, large deviations,\\ and 
integer points on irrational curves}
\author[H. A. {\sc Helfgott}]{{\sc H. A.} HELFGOTT}
\address{H. A. {\sc Helfgott}\\
D\'epartement de math\'ematiques et de statistique\\
Universit\'e de Montr\'eal\\ CP 6128 succ Centre-Ville\\
Montr\'eal, QC H3C 3J7, Canada}
\email{helfgott@dms.umontreal.ca}
\maketitle

\begin{resume}
Soit $f\in \mathbb{Z}\lbrack x\rbrack$ un polyn\^{o}me de degr\'e
$d\geq 3$ sans racines de multiplicit\'e $d$ ou $(d-1)$.
Erd\H{o}s a conjectur\'e que si $f$ satisfait les conditions locales
necessaires alors $f(p)$ est sans facteurs puissances
$(d-1)^{\text{\`emes}}$ pour une infinit\'e de nombres premiers $p$. On prouve 
cela pour
toutes les fonctions $f$ dont l'entropie est assez grande.

On utilise dans la preuve 
un principe de r\'epulsion pour
les points entiers sur les courbes de genre positif et un analogue
arithm\'etique du th\'eor\`eme de Sanov issu de la th\'eorie des
grandes d\'eviations.
\end{resume}

\begin{abstr}
Let $f\in \mathbb{Z}\lbrack x\rbrack$ be a polynomial of degree $d\geq 3$
without roots of multiplicity $d$ or $(d-1)$. 
Erd\H{o}s conjectured that, if $f$ satisfies the necessary local conditions,
then $f(p)$ is free of $(d-1)$th powers for infinitely many primes $p$.
This is proved here for all $f$ with sufficiently high entropy.

The proof serves to demonstrate two innovations: a strong repulsion 
principle for integer points on curves of positive genus, and
a number-theoretical analogue of Sanov's theorem from the theory of
large deviations.
\end{abstr}

\bigskip
\section{Introduction}\label{sec:introdu}
\subsection{Power-free values of $f(p)$, $p$ prime}\label{subs:prima}
Let $f\in \mathbb{Z}\lbrack x\rbrack$ be a polynomial of degree $d\geq 3$
without roots of multiplicity $k$ or greater. It is natural to venture
that there are infinitely many integers $n$ such that $f(n)$ is
free of $k$th powers, unless local conditions fail.
(An integer $a$ is said to be {\em free of $k$th powers} 
if there is no integer $b>1$
such that $b^k|a$.) In fact, such a guess is not only natural, but 
necessary in many applications; for example, we need it to hold with
$k=2$ if we want to approximate the conductor of an elliptic curve in a family
in terms of its discriminant (see \cite{Hell} and \cite{Y}, \S 5, 
for two contexts in which such an approximation is crucial).

Assume an obviously necessary local condition -- namely, 
 that $f(x)\not\equiv 0 \mo p^k$ has a solution in
$\mathbb{Z}/p^k \mathbb{Z}$ for every prime $p$.  
If $k\geq d$, it is easy to prove that there are infinitely many integers
$n$ such that $f(n)$ is free of $k$th powers. If $k<d-1$, proving as much
is a hard and by-and-large open problem. (See \cite{Na}, 
\cite{HN} and \cite{He} for results for $d$ large.) Erd\H{o}s proved that
there are infinitely many $n$ such that $f(n)$ is free of $k$th powers
for $k = d-1$. Furthermore, he conjectured that
there are infinitely
many primes $q$ such that $f(q)$ is free of $(d-1)$th powers, provided that
$f(x)\not\equiv 0 \mo p^{d-1}$ has a solution in
$(\mathbb{Z}/p^{d-1} \mathbb{Z})^*$ for every prime $p$.
This conjecture is needed for applications in which certain
variables are restricted to run over the primes. 
Erd\H{o}s's motivation, however,
may have been the following: there is a difficult diophantine problem implicit
in questions on power-free values -- namely, that 
of estimating the number of integer points on twists of a fixed curve of
positive genus. Erd\H{o}s had managed to avoid this problem
for $k=d-1$ and unrestricted integer argument $n$; if the argument $n$ is
restricted to be a prime $q$, the problem is unavoidable, and
must be solved.

The present paper
proves Erd\H{o}s's conjecture for all $f$ with sufficiently high entropy.
As we will see, even giving a bound of $O(1)$ for the diophantine problem
mentioned above would not be enough; we must mix sharpened diophantine
methods with probabilistic techniques.

We define the
 {\em entropy}\footnote{This is essentially a relative entropy,
appearing as in the theory of large deviations; vd.\ \S \ref{sec:largdev}.} $I_f$ of an irreducible polynomial $f$ over $\mathbb{Q}$
to be
\begin{equation}\label{eq:selbst}
I_f = \frac{1}{|\Gal_f|} \mathop{\sum_{g\in \Gal_f}}_{\lambda_g\ne 0} \lambda_g \log \lambda_g,
\end{equation}
where $\Gal_f$ is the Galois group of the splitting field of $f$
and $\lambda_g$ is the number of roots of $f$ fixed by $g\in \Gal_f$.
(We write $|S|$ for the number of elements of a set $S$.)

\begin{thm}\label{thm:cea}
Let $f\in \mathbb{Z}\lbrack x\rbrack$ be a
polynomial of degree $d$ without roots of multiplicity $\geq k$,
where $k = d-1$ and $d \geq 3$.
If $f$ is irreducible, assume that its entropy $I_f$ is greater than $1$.
Then, for a random\footnote{
Let $S$ be an infinite set of positive integers -- in this case, the primes.
When we say that the probability that a random element $q$ of $S$ satisfy a 
property $P$ is $x$, we mean that the following limit exists and equals $x$:
\[\lim_{N\to \infty} \frac{|\{1\leq q\leq N: \text{$q\in S$
satisfies $P$}\}|}{|\{1\leq q\leq N: q\in S\}|}\]}
 prime $q$, the probability that $f(q)$ be
free of $k$th powers is
\begin{equation}\label{eq:boras}
\prod_{p} \left(1 - \frac{\rho_{f,*}(p^{k})}{p^{k} - 
p^{k-1}}\right),
\end{equation}
where $\rho_{f,*}(p^{k})$ stands for the number of solutions
to $f(x) \equiv 0 \mo p^{k}$ in $(\mathbb{Z}/p^{k} \mathbb{Z})^*$.
\end{thm}
\begin{remark}
The probability (\ref{eq:boras}) is exactly what one would expect from
heuristics: the likelihood that a random prime $q$ be indivisible by
a fixed prime power $p^k$ is precisely $1 - \frac{\rho_{f,*}(p^k)}{p^k -
p^{k-1}}$. The problem is that we will have to work with a set
of prime powers whose size and number depend on $q$.
\end{remark}
It is easy to give a criterion for the non-vanishing of (\ref{eq:boras}).
\begin{cor}\label{cor:orovil}
Let $f\in \mathbb{Z}\lbrack x\rbrack$ be a polynomial of degree $d$
without roots of multiplicity $\geq k$, where $k = d-1$ and $d\geq 3$.
If $f$ is irreducible, assume that its entropy $I_f$ is greater than $1$.
Assume as well that
no $k$th power $m^k$, $m>1$, divides all 
coefficients of $f$, and that
$f(x) \not\equiv 0 \mo p^{k}$ has a solution in
$(\mathbb{Z}/p^{k} \mathbb{Z})^*$ for every $p\leq d+1$. 
Then $f(p)$ is free of $k$th powers for infinitely many primes $p$. 
Indeed, $f(p)$ is free of $k$th powers for a positive proportion of 
all primes.
\end{cor}
\begin{remark}
An irreducible polynomial $f$ of degree $3$, $4$, $5$ or $6$
has entropy greater than $1$ if and only if its Galois group is
one of the following:
\begin{equation}\label{eq:otro}
\begin{aligned}&A_3, C(4), E(4), D(4), C(5),\\ &C(6), D_6(6), D(6), A_4(6), F_{18}(6),
2 A_4(6), F_{18}(6):2,\, F_{36}(6), 2 S_4(6),\end{aligned}\end{equation}
in the nomenclature of \cite{CHM}.
See Table \ref{tab:gd6}. 
Erd\H{o}s's problem
 remains open for 
irreducible polynomials with the following Galois groups:
\[\begin{aligned}
&S_3, A_4, S_4, D(5), F(5), A_5, S_5, \\
&S_4(6d), S_4(6c), L(6), F_{36}(6):2,\,
L(6):2,\, A_6, S_6 .
\end{aligned}\]
\end{remark}
\begin{remark}
We will be able to give bounds on the rate of convergence to
(\ref{eq:boras}): the proportion of primes $q\leq N$ such that
$f(q)$ is free of $k$th powers equals (\ref{eq:boras})$+O((\log N)^{-\gamma})$,
$\gamma>0$. We will compute $\gamma$ explicitly in \S \ref{sec:rater}.
In particular, if $d = 3$ and $\Gal_f = A_3$, then
$\gamma = 0.003567\dotsc$. See \S \ref{sec:rater}, Table \ref{tab:bogo}.
\end{remark}
\begin{remark}
The entropy $I_f$ is greater than $1$ for every normal polynomial $f$
of degree $\geq 3$.
(A polynomial is {\em normal} if one of its roots generates its splitting 
field.) In particular, $I_f>1$ for every $f$ with $\Gal_f$ abelian
and $\deg(f)\geq 3$. We do have $I_f>1$ for many non-normal polynomials
$f$ as well; most of the groups in (\ref{eq:otro}) are Galois groups of
non-normal polynomials. In contrast, for $f$ of degree $d$ with
$\Gal_f = S_d$, the entropy $I_f$ tends to 
$\sum_{k=2}^{\infty} \frac{\log k}{e (k-1)!} = 0.5734028\dotsc$
 as $d\to \infty$. See (\ref{eq:onyo}).
\end{remark}
\begin{remark}
If we can tell whether or not $f, g\in \mathbb{Z}\lbrack x\rbrack$ take
values free of $k$th powers for infinitely many prime arguments, we can
tell the same for $f\cdot g$. In other words, when we work with
a reducible polynomial, the degree and entropy of the
largest irreducible factors of the polynomial matter, rather than
the degree of the polynomial itself. We will take this fact into account
in the statement of the main theorem.
\end{remark}
\begin{table}
\begin{center}
\begin{tabular}{ll|ll|ll}
$\Gal_f$ & $I_f$ 
& $\Gal_f$ & $I_f$ 
& $\Gal_f$ & $I_f$\\ \hline
$A_3$ &   $1.0986123$ 
& $S_3$ & $0.5493061$
& &
\\\hline 
$C(4)$ & $1.3862944$
& $E(4)$ & $1.3862944$ & $D(4)$ & $1.0397208$\\
$A_4$ & $0.4620981$ & $S_4$ & $0.5776227$ & &\\\hline
$C(5)$ & $1.6094379$
& $D(5)$ & $0.8047190$ & $F(5)$ & $0.4023595$\\
$A_5$ & $0.5962179$ & $S_5$ & $0.5727620$ & &\\\hline
$C(6) $ & $1.7917595$
& $D_6(6) $ & $1.7917595$ & $D(6) $ & $1.2424533$\\
$A_4(6) $ & $1.2424533$ & $F_{18}(6) $ & $1.3296613$
&$2A_4(6) $ & $1.3143738$ \\$S_4(6d) $ & $0.9678003$
&$S_4(6c) $ & $0.9678003$ & $F_{18}(6):2 $ & $1.0114043$\\
$F_{36}(6) $ & $1.0114043$ & $2S_4(6) $ & $1.0037605$
&$L(6) $ & $0.5257495$\\
$F_{36}(6):2 $ & $0.9678003$ &$L(6):2 $ & $0.6094484$ & 
$A_6 $ & $0.5693535$\\
$S_6 $ & $0.5734881$ & & & &
\\
\end{tabular}
\caption{Entropies of irreducible polynomials of degree
$3,4,5,6$}\label{tab:gd6}
\end{center}
\end{table}
\subsection{General statement}\label{subs:tojo}
Theorem \ref{thm:cea} holds over many sequences other than the primes.
All we use about the primes is that the proportion of them lying in
a given congruence class can be ascertained, and that they are not much sparser
than a simple sieve majorisation already forces them to be.
\begin{defn}\label{def:boul}
Let $S$ be a set of positive integers. 
We say that $S$ is {\em predictable}
if the limit 
\begin{equation}\label{eq:taois}
\rho(a,m) = \lim_{N\to \infty} \frac{|\{n\in S: n\leq N, n\equiv a \mo m\}|}{|\{n\in S: n\leq N\}|}\end{equation}
exists for all integers $a$, $m>0$.
\end{defn}
The following definition is standard.
\begin{defn}\label{defn:sievdim}
Let $P$ be a set of primes. We say that $P$ is
a {\em sieving set of dimension $\theta$} if
\begin{equation}\label{eq:ontori}
\prod_{w\leq p<z} \left(1 - \frac{1}{p}\right)^{-1} 
\ll \left(\frac{\log z}{\log w}
\right)^{\theta}\end{equation} for all $w$, $z$ with $z>w>1$, where
$\theta\geq 0$ is fixed. 
\end{defn}
We are about to define {\em tight} sets. A tight set is essentially a set 
whose cardinality can be estimated by sieves up to a constant factor.
\begin{defn}\label{defn:oglo}
Let $S$ be a set of positive integers. Let $P$ be a sieving
set with dimension $\theta$.
We say that $S$ is {\em $(P,\theta)$-tight}
 if (a) no element $n$ of $S$ is divisible by any
prime in $P$ smaller than $n^\delta$, where $\delta>0$ is fixed,
 (b) the number of elements of $\{n\in S: n\leq N\}$ is
$\gg N/(\log N)^\theta$ for $X$ sufficiently large.
\end{defn}
In other words, $S$ is a $(P,\theta)$-tight set if the upper bounds on 
its density given by its sieve dimension $\theta$ are tight up to a constant 
factor.
\begin{main}
Let $S$ be a predictable, $(P,\theta)$-tight set.
Let $k\geq 2$. 
Let $f\in \mathbb{Z}\lbrack x\rbrack$ be a
polynomial such that, for every irreducible factor $g$ of $f$,
the degree of $g$ is $\leq k_g + 1$, where $k_g = \lceil k/r_g\rceil$
and $r_g$ is the highest power of $g$ dividing $f$.
Assume that the entropy $I_g$ of $g$ is $> 
(k_g + 1) \theta - k_g$ for every
irreducible factor $g$ of $f$ of
degree exactly $k_g + 1$.

Then, for a random element $q$ of $S$, the probability that $f(q)$ be
free of $k$th powers is
\begin{equation}\label{eq:bormo}
\lim_{z\to \infty} \mathop{\sum_{m\geq 1}}_{p|m \Rightarrow p\leq z}
\mu(m) \mathop{\sum_{0\leq a < m^k}}_{f(a)\equiv 0 \mo m^k}
\rho(a,m^k), 
\end{equation}
where $\rho(a,m)$ is as in (\ref{eq:taois}).
\end{main}
The expression whose limit is taken in (\ref{eq:bormo}) is non-negative and
non-increasing
on $z$, and thus the limit exists.
\begin{example}
The primes are, of course, predictable and $(P,1)$-tight, where $P$
is the set of all primes. Thus, Thm.\ \ref{thm:cea} is a special case of the 
main theorem.
In the general case, if the convergence of (\ref{eq:taois}) is
not too slow, we can obtain bounds for the error term that are of
the same quality
as those we can give in the case of the primes, viz., upper bounds 
equal to the main term times $(1 + O((\log N)^{-\gamma}))$, $\gamma>0$.
\end{example}
\begin{example}
Let $S$ be the set of all sums of two squares. Then
$S$ is predictable and $(P,\frac{1}{2})$-tight, where $P$ is the set of all
primes $p\equiv 3 \mo 4$. Since the entropy (\ref{eq:selbst})
of a polynomial is always positive, we have $I_g > d \cdot 
\frac{1}{2} - (d-1)$ for every irreducible $g$ of degree $d\geq 1$,
and thus we obtain the asymptotic (\ref{eq:bormo}). For this choice of
$S$, the techniques
in \S \ref{sec:trois} and \S \ref{sec:cloud} suffice; the probabilistic
work in \S \ref{sec:largdev} is not needed.
The same is true for
any other $S$ that is $(P,\theta)$ tight with $\theta < \frac{k}{k+1}$.

Note that we are considering sums of two squares counted without
multiplicity. A statement similar to (\ref{eq:bormo}) would 
in fact be true if such
sums are counted with multiplicity; to prove as much is not any harder
than to prove (\ref{eq:bormo}) for $S = \mathbb{Z}$, and thus could be done
with classical sieve techniques. 
\end{example}
\begin{example}
The set of all integers is predictable and $(P,0)$-tight, and thus the
main theorem applies. We will discuss the error terms implicit in
(\ref{eq:bormo}) generally and in detail. Setting $S = \mathbb{Z}$
and $\deg(f)=3$,
we will obtain that the total number of integers $n$ from $1$ to $N$
such that $f(n)$ is square-free equals
$N \prod_p (1 - \rho_f(p)/p)$ plus
$O_f(N (\log N)^{-8/9})$ (if $\Gal_f = A_3$) or
$O_f(N (\log N)^{-7/9})$ (if $\Gal_f = S_3$). 
See Prop.\ \ref{prop:ints}.
The error terms $O_f(N (\log N)^{-8/9})$ and
$O_f(N (\log N)^{-7/9})$ are smaller than those in
 \cite{Hsq}, Thm.\  5.1 (respectively, $O_f(N (\log N)^{-0.8061\dotsc})$ and
$O_f(N (\log N)^{0.6829\dotsc})$), which were, in turn, an improvement over
the bound in \cite{Hoo}, Ch.\ IV
(namely, $O_f(N (\log N)^{1/2})$). Analogous improvements also hold for
square-free values of homogeneous sextic forms; here the strongest
result in the literature so far was \cite{Hsq}, Thm.\  5.2,
preceded by the main theorem in \cite{Gre}.
\end{example}
The main theorem would still hold if the definition of
a $(P,\theta)$-tight set were generalised somewhat.
There is no reason why the sieved-out congruence
class modulo $p$, $p\in P$, should always be the class
$a\equiv 0 \mo p$. One must, however, ensure that, for every
factor $g$ of $f$ with $\deg(g)>1$, we get
$g(a)\not\equiv 0 \mo p$ for all but finitely many of
the sieved-out congruence classes $a \mo p$,
or at any rate for all $p\in P$ outside a set of low
density.
 One may sieve out more than one congruence
class per modulus $p\in P$. 
The number of sieved-out congruence classes
per $p\in P$ need not even be bounded by a constant, but
it ought to be constrained to grow slowly.
\subsection{Plan of attack}
Estimating the number of primes $p$ for which $f(p)$ is {\em not} free
of $k$th powers is the same as estimating the number of
solutions $(t,y,x)$ to $t y^k = f(x)$ with $x$ prime, $x$, $t$, $y$ integers,
$y>1$, and $x$, $t$, $y$ within certain ranges. The solutions
to $t y^k = f(x)$ with $y$ small (or $y$ divisible by a small prime)
can be counted easily. What remains is to bound from above the number
of solutions $(t,y,x)$ to $t y^k = f(x)$ with $x$ and $y$
 prime and $y$ very large --
larger than $x (\log x)^{-\epsilon}$, $\epsilon>0$.
It is intuitively clear (and a consequence
of the $abc$ conjecture; see \cite{Gra}) that such solutions should be
very rare. Bounding them at all non-trivially (and unconditionally) is
a different matter, and the subject of this paper.

{\em Counting integer points on curves.}
Let $C$ be a curve of positive genus $g$. Embed $C$ into its Jacobian $J$.
The abelian group $J(\mathbb{Q})$ is finitely generated; call its rank
$r$.
Map the lattice of rational points of $J$ to $\mathbb{R}^r$ in such a way as 
to send the canonical height to the square of the Euclidean norm. Project
$\mathbb{R}^r\setminus \{0\}$ radially
onto the sphere $S^{r-1}$. Let $P_1$, $P_2$ be
two rational points on $C$ whose difference in $J$ is non-torsion. Mumford's
gap principle amounts in essence to the following statement: if $P_1$ and
$P_2$ are of roughly the same height, then the images of $P_1$ and $P_2$ on
$S^{r-1}$ are separated by an angle of at least $\arccos \frac{1}{g}$.
This separation is not enough for our purposes. We will show that, if $P_1$
and $P_2$ are integral and of roughly the same height, then their images on $S^{r-1}$ are separated by an 
angle of at least $\arccos \frac{1}{2g}$.

The case $g=1$ was already treated in \cite{Hsq}, \S 4.7. The separation
of the points is increased further when, in addition to being 
integral,
 $P_1$ and $P_2$
are near each other in one or more localisations of $C$. This 
phenomenon was already noted in \cite{HV} for $g=1$, as well as in the
case of $P_1$, $P_2$ rational and $g\geq 1$ arbitrary.

In section \S \ref{sec:cloud}, 
we will use the angular separation between integer points
to bound their number. This will be done by means of a lemma
on sphere packings. In our particular
problem, $P_1$ and $P_2$ may generally be taken to be near enough each other in
sufficiently many localisations to bring their separation up to $90^{\circ}
-\epsilon$. We will then have uniform bounds\footnote{Bounds such
as Mumford's $O_{C,L}(\log h_0)$ (\cite{HS}, Thm.\  B.6.5) for the number
of $L$-rational points on $C$ of canonical height up to $h_0$ would
be insufficient: for $C$ fixed and $L$ variable,
 the implied constant is proportional to
$c^{\rnk(J(L))}$, where $c>1$ is a fixed constant. The same is true
of bounds resulting from the explicit version of Faltings' theorem
in \cite{Bo2} -- the bound is then $O_C(7^{\rnk(J(L))})$. Our bound is
$O_C((1+\epsilon)^{\rnk{J(L)}})$ for a typical $t$.
(Here $L = \mathbb{Q}(t^{1/k})$.)} of the 
form $O((\log t)^{\epsilon})$ for the number of points on a typical fibre
$t y^{k} = f(x)$.

{\em Large deviations from the norm.} Let $p$ be a typical prime, i.e.,
a prime outside a set of relative density zero. Suppose that
$t q^{k} = f(p)$ for some prime $q$ and some integer $t < p (\log p)^{\epsilon}$. We can then show that $t$ is, in some ways, a typical integer,
and, in other ways, an atypical one. (We first look at how large
the prime factors of $t$ are, and then at how many there are per
splitting type.)
 The former fact ensures that the
above-mentioned bound $O((\log t)^{\epsilon})$ on the number of points
on $t y^{k} = f(x)$ does hold. The latter fact also works to our advantage:
what is {\em rare} in the sense of being atypical must also be {\em rare}
in the sense of being sparse. (The two senses are one and the same.) Thus
the set of all $t$ to be considered has cardinality much smaller than 
$p (\log p)^{\epsilon}$.

How much smaller? The answer depends on the entropy $I_f$ of $f$. (Hence
the requirement that $I_f>1$ for Theorem \ref{thm:cea} to hold.) Results
on {\em large deviations} measure the unlikelihood
of events far in the tails of probability distributions. We will prove
a variant of a standard theorem (Sanov's; see \cite{Sa} or, e.g.,
\cite{Hol}, \S II.1) where a conditional entropy appears as an exponent.
We will then translate the obtained result into a proposition in number theory,
by means a slight refinement of the Erd\H{o}s-Kac
technique (\cite{EK}). (The refinement is needed because we must translate
the far tails of the distribution, as opposed to the distribution itself.)

Our bounds on the number of $t$'s are good enough when they are better
by a factor of $(\log X)^{\epsilon}$ than the desired bound
of $X/(\log X)$ on the total number of {\em tuples} $(t,q,p)$ satisfying
$t q^k = f(p)$; this is so because our upper bound on the number of points
per $t$ is in general low, viz., $O((\log t)^{\epsilon})$.
\subsection{Relation to previous work}\label{subs:nair}
Using techniques from sieve theory and exponential sums, Hooley 
(\cite{Ho2}, \cite{Ho3})
proved Erd\H{o}s's conjecture for polynomials $f\in \mathbb{Z}\lbrack x\rbrack$
of degree $\deg(f)\geq 51$; for $f$ normal and in a certain sense generic, 
he softened the assumption to $\deg(f)\geq 40$ (\cite{Ho3}, Thms.\ 5, 6).
(The results in the present paper apply to all normal polynomials $f$,
as their entropy is always high enough;
see the comments at the end of \S \ref{sec:pmt}.) Then came a remarkable
advance by Nair \cite{Na}, who, using an approach ultimately derived
from Halberstam and Roth's work on gaps between square-free numbers \cite{HR},
showed that Erd\H{o}s's conjecture holds whenever $\deg(f)\geq 7$. No
other cases of the conjecture have been covered since then.

It is a characteristic common to the rather different approaches in
\cite{Ho2} and \cite{Na} that Erd\H{o}s's conjecture is harder to attack
for $\deg(f)$ small than for $\deg(f)$ large. If one follows the approach
in the present paper, it is not the degree $\deg(f)$ that is crucial,
but the entropy $I_f$: the problem is harder when $I_f$ is small than
when $I_f$ is large.

There are results (\cite{Na}, \cite{HN}, \cite{He}) on values
of $f(n)$ and $f(p)$ free of $k$th powers, where  $k = \deg(f)-2$ or even
lower, provided that $\deg(f)$ be rather high. This is an interesting
situation in which our methods seem to be of no use. 

\subsection{Acknowledgements}
The author would like to thank Alina Cojocaru, E. V. Flynn, Anant Godbole,
Christopher Hall and Anatole Joffe for their patient assistance,
Andrew Granville, for his assistance and encouragement, and Christopher
Hooley, for a dare. 
Thanks are also due to writers of free software (\cite{GAP}) and to an
anonymous referee. 
\section{Notation}
\subsection{Sets}
We denote by $|S|$ the number of elements of a finite set $S$. As is usual,
we say that $|S|$ is the {\em cardinality} of $S$.
\subsection{Primes}
By $p$ (or $q$, or $q_1$, or $q_2$) we shall always mean a prime.
We write $\omega(n)$ for the number of prime divisors of an integer
$n$, and $\pi(N)$ for the number of primes from $1$ up to $N$. Given two
integers $a$, $b$, we write $a|b^{\infty}$ if all prime divisors of $a$
also divide $b$, and $a\nmid b^{\infty}$ if there is some prime divisor of $a$
that does not divide $b$. We define $\gcd(a,b^{\infty})$ to be the largest
positive integer divisor of $a$ all of whose prime factors divide $b$.
\subsection{Number fields}
Let $K$ be a number field. We write $\overline{K}$ for an
algebraic closure of $K$. Let $M_K$ be the set of places of $K$.
We denote the completion of $K$
at a place $v\in M_K$ by $K_v$. 
If $f\in \mathbb{Q}\lbrack x\rbrack$ is an irreducible polynomial, let
$\Gal_f$ be the Galois group of the splitting field of $f$.

If $\mathfrak{p}$
 is a prime ideal of $K$, 
we denote the place corresponding to $\mathfrak{p}$ by
$v_{\mathfrak{p}}$. Given $x\in K^*$, we define $v_{\mathfrak{p}}(x)$
to be the largest integer $n$ such that $x\in \mathfrak{p}^n$.
Define absolute values $|\cdot |_{v_p}$ on $\mathbb{Q}$ by
$|x|_{v_p} = p^{-v_p(x)}$. If $w$ is a place of $K$, and $v_p$ is the
place of $\mathbb{Q}$ under it, then $|\cdot |_w$ is normalised
so that it equals $|\cdot |_{v_p}$ when restricted to $\mathbb{Q}$.

Given a positive integer $n$ and a conjugacy class $\langle g\rangle$
in $\Gal(K/\mathbb{Q})$, we write
$\omega_{\langle g\rangle}(n)$ for 
$\sum_{p|n,\,\text{$p$ unramified},\, \Frob_p=\langle g\rangle} 1$,
where $\Frob_p$ denotes the Frobenius element of 
$p$ in $K/\mathbb{Q}$. 
\subsection{Curves}
As is usual, we denote local
heights with respect to a divisor $D$ by $\lambda_{D,v}$,
and the global height by $h_{D}$. Let $C$ be a curve over a local field
$K_w$, and let $R$ be a point on $C$. We then say that a point $P$ on $C$ is 
{\em integral
with respect to $(R)$} if $f(P)$ is in the integer ring of $K_w$
 for every rational function $f$ on $C$ without poles outside $R$.
Given a curve $C$ over a number field $K$, a set of places $S$ including
all archimedean places, and a point $R$ on $C$, we say that a point
$P$ on $C$ is {\em $S$-integral with respect to $(R)$} if
$P$ is integral on $C\otimes K_w$ with respect to $(R)$ for every place
$w\notin M_K\setminus S$.

\subsection{Functions} We will write $\exp(x)$ for $e^x$. We define
$\li(N) = \int_{2}^N \frac{d x}{\ln x}$.
\subsection{Probabilities} We denote by $\Prob(E)$ the
probability that an event $E$ takes place.
\section{Repulsion among integer points on curves}\label{sec:trois}
Consider a complete non-singular curve $C$ of genus 
$g\geq 1$ over a number field $K$. Embed $C$ in its Jacobian $J$
by means of the map $P\mapsto \Cl(P) - (P_0)$, where $P_0$
is a fixed arbitrary point on $C$. Let $\langle \cdot, \cdot \rangle
: J(\overline{K}) \times J(\overline{K}) \to \mathbb{R}$,
$| \cdot | : J(\overline{K}) \to \mathbb{R}$ be the inner product
and norm induced by the canonical height corresponding to the
theta divisor $\theta\in \Div(J)$. Denote by $\Delta$ the
diagonal divisor on $C\times C$. 

\begin{thm}\label{thm:rheinsilber}
Let $K$ be a number field.
We are given a complete non-singular curve $C/K$
of genus $g\geq 1$ with an embedding $P\mapsto \Cl(P) - (P_0)$
into its Jacobian $J(C)$.
Let $R$ be a point on $C$, and let 
$S$ be any set of places of $K$ including all archimedean places.
Let $L/K$ be an extension of degree $d$; write $S_L$
for the sets of places of $L$ above $S$.

Then, for any two distinct points $P, Q \in C(L)$
that are $S_L$-integral with
respect to $(R)$,
\begin{equation}\label{eq:gharab}\begin{aligned}
\langle P, Q\rangle &\leq \frac{1 + \epsilon}{2 g} (|P|^2 + |Q|^2) -
\frac{1 - \epsilon}{2 g} \max(|P|^2,|Q|^2) \\ &+ \frac{1}{2} \delta
- \frac{1}{2} \sum_{w\in M_L\setminus S_L} d_w
\max(\lambda_{\Delta,w}(P,Q),0) 
+ O_{C,K,\epsilon,d,R,P_0}(1)\end{aligned}
\end{equation}
for every $\epsilon>0$, where
\begin{equation}\label{eq:eat}\delta = 
\sum_{w\in S_L} d_w (
\max(\lambda_{(R),w}(P),\lambda_{(R),w}(Q)) -
\min(\lambda_{(R),w}(P),\lambda_{(R),w}(Q)))\end{equation}
and $d_w = \lbrack L_w : \mathbb{Q}_p\rbrack/\lbrack L : \mathbb{Q}\rbrack$,
where $p$ is the rational prime lying under $w$.
\end{thm}
The fact that the error term
$O_{C,K,\epsilon,d,R,P_0}(1)$ does not depend on $L$ will be crucial to
our purposes.
\begin{proof}
We may state Mumford's gap principle as follows:
\begin{equation}\label{eq:ajli}
2 g \langle P, Q\rangle \leq (1 + \epsilon) (|P|^2 + |Q|^2) - 
g h_{\Delta}(P,Q) + O_{C,P_0,\epsilon}(1).\end{equation}
(See, e.g., \cite{La}, Thm.\  5.11, or \cite{HS}, Prop.\ 
B.6.6\footnote{There is a 
factor of $\frac{1}{2}$ missing before $h_{C\times C,
\Delta}(P,Q)$ in \cite{HS}; cf.\  \cite{HS}, top of p.\ 218. 
Note that, as \cite{HS} states, (\ref{eq:ajli}) 
is valid even for $g=1$.}.) Our task is to show that
the contribution of $g h_{\Delta}(P,Q)$ must be large. Without it,
we would have only the angle of $\arccos \frac{1}{2 g}$ mentioned
in the introduction, as opposed to an angle of $\arccos \frac{1}{g}$.
(We would not, in fact, be able to do any better than
$\arccos \frac{1}{2 g}$ if we did not know that
$P$ and $Q$ are integral.)

We will argue that, since $P$ and $Q$ are $S$-integral,
their heights are made almost entirely out of the contributions
of the local heights
$\lambda_v$, $v\in S$, and that these contributions, minus $\delta$,
are also present in $h_{\Delta}(P,Q)$. Then we will examine the contribution
of the places outside $S$ to $h_{\Delta}(P,Q)$; the expression
 $\sum_{w\in M_L\setminus S_L} 
d_w \max(\lambda_{\Delta,w}(P,Q),0)$ will give a lower bound to this
contribution.

Write $h_{\Delta}(P,Q) = \sum_w d_w \lambda_{\Delta,w}(P,Q) +
O_C(1)$ (as in, say, \cite{HS}, Thm.\ B.8.1(e)).
By \cite{Si},
Prop.\ 3.1(b), every $\lambda_{\Delta,w}$ satisfies
\begin{equation}\label{eq:aghur}
\lambda_{\Delta,w}(P,Q) \geq \min(\lambda_{\Delta,w}(R,P),
\lambda_{\Delta,w}(R,Q)) .\end{equation}
We have 
\begin{equation}\label{eq:anyar}
\lambda_{\Delta,w}(R,P) = \lambda_{(R),w}(P),\;\;\;
\lambda_{\Delta,w}(R,Q) = \lambda_{(R),w}(Q)\end{equation}
 by \cite{Si}, Prop.\ 3.1(d). Thus $h_{\Delta}(P,Q)$ is at least
\begin{equation}\label{eq:twiddle}
\max\left(\sum_{w\in S_L} d_w \lambda_{(R),w}(P),
\sum_{w\in S_L} d_w \lambda_{(R),w}(Q)\right)
- \delta + \sum_{w\in M_L\setminus S_L} d_w \lambda_{\Delta,w}(P,Q) 
\end{equation}
plus $O_C(1)$.

We must first show that $\sum_{w\in S_L} d_w \lambda_{(R),w}(P)$
equals $h_{(R)}(P)$ plus a constant, and similarly for
$h_{(R)}(Q)$. 
Let $w\in M_L\setminus S_L$.
If $w$ is non-archimedean and $C$ has good reduction at $w$, 
the height $\lambda_{(R),w}(P)$ (resp.x $\lambda_{(R),w}(Q)$)
is given by the intersection product $(R \cdot P)$ 
(resp. $(R \cdot Q)$)
on the reduced curve $C \otimes \mathbb{F}_w$ 
(\cite{Gro}, (3.7)).
Since $P$ and $Q$ are integral with respect to $(R)$, both $(R \cdot P)$
and $(R \cdot Q)$ are $0$. Hence
\[\lambda_{(R),w}(P) = \lambda_{(R),w}(Q)=0 .\]
Consider now the case where $w$ is archimedean or
$C$ has bad reduction at $w$. Choose any rational function
$f$ on $C$ whose zero divisor is a non-zero multiple of $R$. 
Since $P$ and $Q$ are integral, both
 $|f(P)|_w$ and $|f(Q)|_w$ are $\geq 1$.
 By functoriality (\cite{HS},
Thm.\  B.8.1(c)) and the fact that, under the standard definition
of the local height on the projective line, $\lambda_{(0),w}(x) = 0$ for
any $x=(x_0,x_1)$ on $\mathbb{P}^1$ with $\left|\frac{x_0}{x_1}\right|_w\geq 1$
(see, e.g., \cite{HS}, Ex.\ B.8.4), 
 it follows that
\begin{equation}\label{eq:rece}\lambda_{(R),w}(P) = O_{C,R,L_w}(1),\;\;\; 
\lambda_{(R),w}(Q) = O_{C,R,L_w}(1) .\end{equation}
Every place $w$ of $L$ that is archimedean or 
of bad reduction must lie above a place $v$ of $K$ that is
archimedean or of
bad reduction. Since there are only finitely many
such $v$, and finitely many extensions $w$ of degree at most $d$ of
each of them (see, e.g., \cite{LAnt}, Ch.\ II, Prop.\ 14), we conclude that
\begin{equation}\label{eq:ojul}\begin{aligned}
h_{(R)}(P) &= \sum_{w\in S_L}
d_w \lambda_{(R),w}(P) + O_{C,R,K,d}(1),\\
h_{(R)}(Q) &= \sum_{w\in S_L}
d_w \lambda_{(R),w}(Q) + O_{C,R,K,d}(1) .\end{aligned}\end{equation}

Now, again by an expression in terms of intersection products, 
$\lambda_{\Delta,w}(P,Q)$ is non-negative
at all non-archimedean places $w$ where $C$ has good reduction,
and, by (\ref{eq:aghur}), (\ref{eq:anyar}) and (\ref{eq:rece}),
it is bounded below by $O_{C,\Delta,L_w}(1)$ at all other places $w$.
We use both these facts and (\ref{eq:ojul}) to bound
(\ref{eq:twiddle}) from below, and we obtain that $h_{\Delta}(P,Q)$
is at least
\[\begin{aligned}
\max(h_{(R)}(P),h_{(R)}(Q)) &- \delta +
\sum_{w\in M_L\setminus S_L} d_w
\max(\lambda_{\Delta,w}(P,Q),0) \\ &+ 
 O_{C,K,R,d}(1) .\end{aligned}\]
By the argument at the bottom of p.\ 217 in \cite{HS}
with $R$ instead of $P_0$,
we have
 $|P|^2 \leq g (1 + \epsilon) h_{(R)}(P) + O_{C}(1)$, 
$|Q|^2 \leq g (1 + \epsilon) h_{(R)}(Q) + O_{C}(1)$.
We apply (\ref{eq:ajli}) and are done.
\end{proof}
The general applicability of Thm.\
\ref{thm:rheinsilber} is somewhat limited by the presence of
a term $O_{C,K,\epsilon,d,R,P_0}(1)$ depending on the curve $C$.
(For the application in this paper, it will be good enough to know that  
$O_{C,K,\epsilon,d,R,P_0}(1)$ does not depend on $L$, but just on its
degree $d = \deg(L/K)$.)
The main obstacle to a uniformisation in the style of \cite{HV},
Prop.\ 3.4, seems to be a technical one: we would need
 explicit expressions for local heights at places of bad reduction,
and the expressions available for genus $g>1$ are not
explicit enough.

\section{Counting points on curves}\label{sec:cloud}
We must now clothe \S \ref{sec:trois} in concrete language for the
sake of our
particular application. Since the field $L$ in Thm.\ 
\ref{thm:rheinsilber} will now be of the special form
 $L = \mathbb{Q}(t^{1/k})$, we will be able to give a 
 bound $\rnk(J(L))$ in terms of the number of prime divisors of $t$
by means of a simple descent argument. We will
then combine Thm.\ \ref{thm:rheinsilber} with sphere-packing results
to give a low bound
((\ref{eq:louise})) on
the number of solutions to $t y^{d-1} = f(x)$
with $t$ fixed and typical.

\begin{lem}\label{lem:sphpack}
Let $A(n,\theta)$ be the maximal number of points that can be arranged on
the unit sphere of
$\mathbb{R}^n$ with angular separation no smaller than $\theta$. Then, for
$\epsilon>0$, 
\[ \lim_{n\to \infty} 
\frac{1}{n} \log_2 A(n,\frac{\pi}{2} - \epsilon) = O(\epsilon) .\]
\end{lem}
\begin{proof}
Immediate from standard
 sphere-packing bounds; see \cite{KL} (or the expositions in \cite{Le}
and \cite{CS}, Ch. 9) for stronger statements. In particular, $O(\epsilon)$
could be replaced by $O(\epsilon^2 \log \epsilon^{-1})$.
\end{proof}

When we speak of the rank of a curve over a field $K$, 
we mean, as is usual, the rank 
of the abelian group of $K$-rational points on its Jacobian.
\begin{lem}\label{lem:angr}
Let $f\in \mathbb{Z}\lbrack x\rbrack$ be a polynomial
of degree $d\geq 3$ without repeated roots. 
Let $p$ be a prime that does not divide $d$. Let $K/\mathbb{Q}$
be a number field.
 Then, for any non-zero integer $t$, the 
curve
\[C_t : t y^p = f(x)\]
has rank over $K$ at most $d (p-1) \lbrack K : \mathbb{Q}\rbrack
 \cdot \omega(t) + O_{K,f,p}(1)$.
\end{lem}
\begin{proof}
Let $J$ be the Jacobian of $C_t$.
Let $\phi$ be the endomorphism $1-\tau$ of $J$,
where $\tau$ is the map on $J$ induced by the map
$(x,y) \mapsto (x,\zeta_p y)$ on $C_t$.
By  \cite{Sch}, Cor.\ 3.7 and Prop.\ 3.8,
\[\rnk_{\mathbb{Z}}(J(K))\leq \frac{p-1}{\lbrack K(\zeta_p) :
K\rbrack}
\rnk_{\mathbb{Z}/p \mathbb{Z}}(J(K(\zeta_p))/\phi 
J(K(\zeta_p))) .\]
By the proof of the weak Mordell-Weil theorem, 
$J(K(\zeta_p))/\phi J(K(\zeta_p))$
injects into $H^1(K(\zeta_p),J\lbrack \phi\rbrack; S)$,
where $S$ is any set of places of $K(\zeta_p)$
containing all places where $C_t$ has bad reduction in addition to
a fixed set of places.
By \cite{Sch}, Prop.\ 3.4, the rank of
$H^1(K(\zeta_p),J\lbrack \phi\rbrack; S)$ 
over $\mathbb{Z}/p \mathbb{Z}$ is no greater
than the rank of $L(S_L,p)$, where $L = K(\zeta_p)\lbrack T
\rbrack /(t^{p-1} f(T))$ and $S_L$ is the set of places of $L$
lying over $S$. (Here $L(S_L,p)$ is the subgroup of $L^*/L^{* p}$
consisting of the classes $\mo L^{* p}$ represented by 
elements of $L^*$ whose valuations at all places outside
$S_L$ are trivial.)
As the roots of $t^{p-1} f(x) = 0$ are independent of $t$, so is $L$.
Thus, the rank of $L(S_L,p)$ is 
$|S_L| + O_{K,f,p}(1) \leq d \cdot |S| + O_{K,f,p}(1)$, 
where the term $O_{K,f,p}(1)$ comes from the size of the
class group of $L$ and from the rank of the group of units of $L$.
The number of places of bad reduction of $C_t$ over
$K(\zeta_p)$ is at most $\lbrack K(\zeta_p) : \mathbb{Q}\rbrack 
\omega(t) + O_{K,f,p}(1)$,
where $O_{K,f,p}(1)$ stands for the number of prime ideals of
$K(\zeta_p)$ dividing the discriminant of $f$. 
The statement follows.
\end{proof}
\begin{prop}\label{prop:notung}
Let $f\in \mathbb{Z}\lbrack x\rbrack$ be a 
polynomial of degree $d\geq 3$ with no repeated roots.
Let $k\geq 2$ be an integer such that $k\nmid d^{\infty}$.
 Let $t\leq X$ be a positive integer. Suppose that $t$
has an integer divisor $t_0\geq X^{1-\epsilon}$, $\epsilon>0$,
such that $\gcd(t_0,(k (\Disc f))^{\infty})$ is less than a constant $c$.
Then the number of integer solutions to $t y^k = f(x)$ with 
$X^{1 - \epsilon} < x\leq X$ is at most 
\begin{equation}\label{eq:louise}
O_{f,k,c,\epsilon}\left(
 e^{O_{f,k}(\epsilon \omega(t))} \prod_{p|t_0} \rho(p)
\right),
\end{equation}
where $\omega(t)$ is the number of prime divisors of $t$ and
$\rho(p)$ is the number of solutions to $f(x) \equiv 0 \mo p$.
\end{prop}
The divisor $t_0|t$ here plays essentially
the same role as the ideal $\mathscr{I}$
in the proof of Thm.\ 3.8 in \cite{HV}. The main difference is that, in
our present case, the congruence $f(x)\equiv 0 \mo t_0$ makes the
cost of considering all possible congruence classes $x \mo t_0$
quite negligible.

The case $k|d^{\infty}$, $k$ not a power of $2$ 
(or, in general, $k$ such that $\gcd(k,d)>2$) is covered
by the recent work of Corvaja and Zannier (\cite{CZ}, Cor.\ 2). Be that
as it may, we will need only the case $k\nmid d^{\infty}$, and thus will not
use \cite{CZ}. 
We could, at any rate, modify Lem.\ \ref{lem:angr} 
to cover the case $p|d$ by using \cite{PSc}, \S 13, instead of \cite{Sch},
\S 3.
Proposition \ref{prop:notung} would then cover the case $k|d^{\infty}$.

\begin{proof}[Proof of Prop.\ \ref{prop:notung}]
Choose a prime $q$ dividing $k$ but not $d$. Define $K = \mathbb{Q}$,
$L = \mathbb{Q}(t^{1/q})$. Let $S$ and $S_L$
be the sets of archimedean places of $\mathbb{Q}$ and $L$, respectively.
Consider the curve $C:y^k = f(x)$. Denote the point at infinity on $C$ by
$\infty$. Embed $C$ into
its Jacobian by means of the map $P\mapsto (P) - (\infty)$.

Now consider any two distinct solutions $(x_0,y_0)$, $(x_1,y_1)$
to $t y^q = f(x)$ with $X^{1-\epsilon} \leq x_0,x_1\leq X$ and
$x_0\equiv x_1 \mo t_0$. Then the points $P=(x_0,t^{1/q} y_0)$,
$Q=(x_1,t^{1/q} y_1)$ on $C$ are integral with respect to 
$S_L$ and $(\infty)$. We intend to apply Thm.\ \ref{thm:rheinsilber},
and thus must estimate the quantities on the right side of (\ref{eq:gharab}).

By the additivity and functoriality of
the local height (\cite{HS}, Thm B.8.1, (b) and (c)) and the fact that
the point at infinity on $\mathbb{P}^1$ lifts back to $q\cdot
\infty$ on $C$ under the map $(x,y)\mapsto x$,
\[\begin{aligned}
\frac{1 - \epsilon}{q} \log X + O_{f,q,w}(1)
&\leq \lambda_{\infty,w}(P) \leq
\frac{1 + \epsilon}{q} \log X + O_{f,q,w}(1),\\
\frac{1 - \epsilon}{q} \log X + O_{f,q,w}(1)
 &\leq \lambda_{\infty,w}(Q) \leq
\frac{1 + \epsilon}{q} \log X + O_{f,q,w}(1)
\end{aligned}\]
for $w\in S_L$. We know that
$||P|^2 - g h_{\infty}(P)| \leq \epsilon h_{\infty}(P) + O_C(1)$
and $|Q|^2 - g h_{\infty}(Q)|\leq \epsilon h_{\infty}(Q) + O_C(1)$
(vd., e.g., the argument at the bottom of p.\ 217 in \cite{HS}). Hence
\begin{equation}\label{eq:illo}\begin{aligned}
\frac{(1 - \epsilon)^2}{q} g \log X + O_{f,q}(1)&\leq |P|^2\leq
\frac{(1 + \epsilon)^2}{q} g \log X + O_{f,q}(1),\\
\frac{(1 - \epsilon)^2}{q} g \log X + O_{f,q}(1) &\leq |Q|^2\leq
\frac{(1 + \epsilon)^2}{q} g \log X + O_{f,q}(1)\end{aligned}\end{equation}
Since $P = (x_0, t^{1/q} y_0)$, $Q = (x_1, t^{1/q} y_1)$ and $t_0|t$,
we have that,
for every non-archimedean place $w\in M_L$ where $C$ has good reduction,
$\lambda_{\Delta,w}(P,Q) \geq - \log |t_0^{1/q}|_w$ (see, e.g.,
\cite{La3}, p.\ 209). Thus, for every prime
$p|t_0$ where $C$ has good reduction,
\[\sum_{w|p} d_w \lambda_{\Delta,w}(P,Q) \geq - \sum_{w|p} d_w
\log |t_0^{1/q}|_w = \frac{1}{q} p^{v_p(t_0)} ,\]
where $d_w = \lbrack L_w:\mathbb{Q}_p\rbrack/\lbrack L : \mathbb{Q}\rbrack$.
We apply Thm.\ \ref{thm:rheinsilber} and obtain
\begin{equation}\label{eq:armad}\begin{aligned}
\langle P,Q\rangle &\leq 
\frac{(1 + \epsilon)^3}{2 g q} (g \log X + g \log X) \\ &- 
\frac{(1 - \epsilon)^3}{2 g q} g
\log X + \frac{\epsilon}{q} \log X - \frac{1}{2 q} 
\sum_{p} \log p^{v_p(t_0)} + O_{f,k,\epsilon}(1) 
\\ &= O\left( \frac{\epsilon}{q} \log X\right) + O_{c,f,k,\epsilon}(1)
,\end{aligned}\end{equation}
where we use the facts that $t_0\geq X^{1-\epsilon}$ and that the
sum of $\log p^{v_p(t_0)}$ over all primes $p$ of bad reduction is bounded 
above by the constant $c$.

By (\ref{eq:illo}) and (\ref{eq:armad}), we conclude that,
for $X$ large enough (in terms of $c$, $f$, $k$ and $\epsilon$), 
$P$ and $Q$ are separated by
an angle of at least $\pi/2 - O_{f,k}(\epsilon)$ 
in the Mordell-Weil lattice 
$J(L)$ endowed with the inner product $\langle \cdot, \cdot\rangle$ induced by the theta divisor.
 By Lemma \ref{lem:sphpack}, there
can be at most $e^{O(\epsilon r)}$ 
points in $\mathbb{R}^r$
separated by angles of at least $\pi/2 - O(\epsilon)$. Since the rank
$r$ of $J(L)$ is bounded from above by $O_{f,k}(\omega(t))$
(Lemma \ref{lem:angr}), it follows that there can be at most
$e^{O_{f,k}(\epsilon \omega(t))}$ 
points placed as $P$ and $Q$
are, viz., satisfying $X^{1-\epsilon}\leq x\leq X$ and
having $x$-coordinates congruent to each other modulo $t_0$.
Since $t y^q = f(x)$ implies $f(x) \equiv 0 \mo t_0$, there are at most
$O_f(\prod_{p|t_0} \rho(p))$ congruence classes 
modulo $t_0$ into which $x$ may fall.
\end{proof}
\section{The probability of large deviations}\label{sec:largdev}
Our task in this section will be to translate into number theory 
a statement (Sanov's theorem, \cite{Sa}) on the probability of unlikely
events. (If a die is thrown into the air $n$ times,
where $n$ is large, what is the order of the probability
that there will be fewer than $\frac{n}{10}$ ones
and more than $\frac{n}{5}$ sixes? The central limit theorem does
not yield the answer; it only tells us that the probability goes to
zero as $n$ goes to infinity.)
The translation resembles the argument in \cite{EK}, though some of
the intermediate results must be sharpened.

Let $J$ be a finite index set. For 
$\vec{c}, \vec{x} \in (\mathbb{R}_0^+)^J$,
define
\begin{equation}\label{eq:tuastort}
B_{\vec{c},\vec{x}} = \{\vec{y} \in (\mathbb{R}_0^+)^J :
\sgn(y_j-x_j) = \sgn(x_j-c_j)\;\;\;\; \text{$\forall j\in J$ s.t.
$x_j\ne c_j$} \} ,\end{equation}
where $\sgn(t)$ is as follows:
 $\sgn(t)=1$ if $t>0$, $\sgn(t)=-1$ if $t<0$, and $\sgn(t)=0$ if $t=0$.
In other words, $B_{\vec{c},\vec{x}}$ is the set of all vectors $\vec{y}$ that
are no closer to $\vec{c}$ 
than $\vec{x}$ is: $y_j<x_j$ if $x_j<c_j$, and
$y_j> x_j$ if $x_j> c_j$. We also define
\begin{equation}\label{eq:jacbo}
I_{\vec{c}}(\vec{x}) = 1 - \sum_{j\in J} x_j + \sum_{j\in J} x_j 
 \log \frac{x_j}{c_j} .\end{equation}
We adopt the convention that, if $c_j = 0$, then
 $\log \frac{x_j}{c_j} = \infty$, unless $x_j$ also equals $0$, in which
case we leave $\log \frac{x_j}{c_j}$ undetermined and
take $x_j \log \frac{x_j}{c_j}$ to be $0$.

The following is a variant of Sanov's theorem.
\begin{prop}\label{prop:san}
Let the rational primes be partitioned into
$\{P_j\}_{j\in J}$, $J$ finite,
so that, for every $j\in J$, we have the asymptoptic
$\sum_{p\in P_j,\, p\leq N} 1/p \sim r_j \log \log N$,
where $\vec{r} \in (\mathbb{R}_0^+)^d$. 
Let $\{X_p\}_{\text{$p$ prime}}$ be jointly
independent random variables with
values in $(\mathbb{R}_0^+)^d$ defined by
\begin{equation}\label{eq:carne}
X_p = \begin{cases} e_j & \text{with probability $s_j/p$,}\\
0 &\text{with probability $1 - s_j/p$,}\end{cases}\end{equation}
where $\vec{s} \in (\mathbb{R}_0^+)^d$,
$e_j$ is the $j$th unit vector in $\mathbb{R}^J$ and
$j\in J$ is the index such that $p\in P_j$. 

Define $\vec{c} \in (\mathbb{R}_0^+)^J$ by
$c_j = r_j s_j$. 
Then, for all $\vec{x} \in (\mathbb{R}_0^+)^J$, 
\[\lim_{n\to \infty} \frac{1}{\log \log n}
\log \Prob\left(\frac{1}{\log \log n} \sum_{p\leq n} \delta_{X_p}
\in B_{\vec{c},\vec{x}}\right)
= - I_{\vec{c}}(\vec{x}),\]
where $I_{\vec{c}}(\vec{x})$ is as in (\ref{eq:jacbo})
and $\delta_{\vec{x}}$ denotes the point mass at $\vec{x}\in \mathbb{R}^d$.
\end{prop}
\begin{proof}
For $m>0$, let $Z_m = \frac{1}{m} \sum_{p\leq e^{e^m}} \delta_{X_p}$.
Define $\phi_m(\vec{\,t}\,) = \mathbb{E}\left(e^{\langle \vec{t}\,,Z_m\rangle}\right)$
for $\vec{t} \in \mathbb{R}^J$. 
Then
\[\phi_m(m \vec{t}\,) = \mathbb{E}\left(e^{\langle m \vec{t}, Z_m\rangle}\right)
= \prod_{j\in J} \mathop{\prod_{p\leq e^{e^m}}}_{p\in
P_j} \left(\left(1 - \frac{s_j}{p}\right) + \frac{s_j}{p} e^{t_j}\right)
.\]
Define $\Lambda(\vec{t}\,) = \lim_{m\to\infty} \frac{1}{m} \log 
\phi_m(m \vec{t}\,)$. We
obtain
\[\Lambda(\vec{\,t}\,) = \sum_{j\in J}
\lim_{m\to \infty} \frac{1}{m} \mathop{\sum_{p\leq e^{e^m}}}_{p\in
P_j} \log \left( 1 + \frac{s_j}{p} (e^{t_j} - 1)\right)
= \sum_{j\in J} c_j (e^{t_j} - 1) .\]
Write $\Lambda^*(\vec{y})$ for the Legendre transform 
$\sup_{\vec{t}\in \mathbb{R}^J} (\langle \vec{y},\vec{\,t}\,\rangle -
\Lambda(\vec{\,t}\,) )$ of $\Lambda(\vec{t}\,)$. 
For $\vec{y}\in (\mathbb{R}_0^+)^J$ with $y_j=0$ for every
$j\in J$ with $c_j=0$,
the maximum
of $\langle \vec{y},\vec{t}\,\rangle - \Lambda(\vec{\,t}\,)$ is attained 
at all $\vec{t} \in (\mathbb{R}_0^+)^J$ such that
$t_j = \log \frac{y_j}{c_j}$ for every $j\in J$ with $c_j\ne 0$. Thus, 
$\inf_{\vec{y} \in B_{\vec{c},\vec{x}}} \Lambda^*(\vec{y})$ equals
\[\mathop{\inf_{\vec{y} \in B_{\vec{c},\vec{x}}}}_{c_j=0 \Rightarrow y_j=0}
 \left(1 - \mathop{\sum_{j\in J}}_{c_j\ne 0} y_j + 
\mathop{\sum_{j\in J}}_{c_j\ne 0} y_j
\log \frac{y_j}{c_j}   \right) = 1 - \sum_{j\in J} x_j + \sum_{j\in J} x_j 
\log 
\frac{x_j}{c_j} = I_{\vec{c}}(\vec{x}) .\]
(The equation is valid even if $c_j=0$ for some $j\in J$, thanks
to our convention that $x_j \log(x_j/c_j) = 0$ when $x_j=c_j=0$.
For $\vec{y} \in (\mathbb{R}_0^+)^J$ such that $y_j\ne 0$, $c_j=0$ for some
$j\in J$, the function $\vec{t} \mapsto \langle \vec{y},\vec{t}\rangle - 
\Lambda(\vec{t})$ is unbounded above, and so $\Lambda^*(\vec{y}) =\infty$.)
By the G\"artner-Ellis theorem (see, e.g., \cite{Hol}, Thm.\  V.6, or
\cite{DZ},  Thm.\  2.3.6), we conclude that
\[\lim_{m\to \infty} \frac{1}{m} \log(\Prob(Z_m \in B_{\vec{c},\vec{x}})) = - 
I_{\vec{c}}(\vec{x}) .\]
\end{proof}
The following lemma serves a double purpose. It is a crucial step in the 
translation of a probabilistic large-deviation result (in our case,
Prop.\ \ref{prop:san}) into arithmetic (cf.\ \cite{EK}, Lemma 4). Later, it will
also allow us to apply Prop.\ \ref{prop:notung} in such as way as to get
a bound of $(\log d)^{\epsilon}$ for the number of integral points of moderate
height on the curve $d y^{r-1} = f(x)$, where $d$ is any integer 
outside a sparse exceptional set.
\begin{lem}\label{lem:spex}
Let $f\in \mathbb{Z}\lbrack x\rbrack$ be a polynomial. Then, for any
$A>0$, $\epsilon>0$, there is a function 
$\delta_{f,A,\epsilon}:(e,\infty)\to \lbrack 0,1\rbrack$ with
$|\log \delta(x)| < \epsilon \log \log x$ and
$\delta(x) = o(1/\log \log x)$, such that, for all but
$O_{f,A,\epsilon}(N (\log N)^{-A})$ integers $n$ between $1$ and $N$,
\begin{enumerate}
\item\label{it:oyop} $\prod_{p|f(n): p\leq N^{\delta(N)}} p < N^{\epsilon}$,
\item\label{it:egg}
 $\sum_{p|f(n): p> N^{\delta(N)}} 1 +\sum_{p^2|f(n): p\leq N^{\delta(N)}} 1\;< \epsilon \log \log N$.
\end{enumerate}
\end{lem}
In other words, the bulk in number of the divisors is on one side,
and the bulk in size is on the other side. All but very few
of the prime 
divisors of a typical number 
are small, but their product usually amounts to very little. 
\begin{proof}
Define $\gamma(n) = \prod_{p|f(n): p\leq N^{\delta(N)}} p$.
Let $\delta(x) = (\log x)^{-\epsilon/ r c^{2 r}}$, where
$r = \deg(f)$ and $c$ will be set later in terms of $A$ and $\epsilon$. 
Then, for any positive integer
$k$ and all $N$ such that $\delta(N) < \frac{1}{k}$,
\[\begin{aligned}
\sum_{1\leq n\leq N} (\log \gamma_N(n))^k &= \sum_{1\leq n\leq N}
\left(\sum_{p|f(n): p\leq N^{\delta(N)}} \log p \right)^k \\
&\ll_{k,f} 
N \max_{1\leq j \leq k} \left(\sum_{p\leq N^{\delta(N)}}
\frac{r \log p}{p}\right)^j \\
&\ll_{r,k} N (\log N^{\delta(N)})^k = N (\log N)^{k - \epsilon k/r c^{2 r}} .
\end{aligned}\]
Setting $k = \lceil A r c^{2 r}/\epsilon\rceil$, we obtain that there are
$O_{c,f,A,\epsilon}(N (\log N)^{-A})$ integers $n$ from $1$ to $N$ such that
$\gamma_N(n)\geq N^{\epsilon}$. Thus (\ref{it:oyop}) is fulfilled.

Clearly \[c^{\omega(f(n)/\gamma_N(n))} =
(c^{2 r})^{\omega(f(n)/\gamma_N(n))/(2 r)} \leq
\mathop{\max_{\text{$d$ sq.-free},\; d\leq C \sqrt{N}}}_{d|f(n)/\gamma_N(n)}
c^{2 r \omega(d)},\]
where $C$ is the absolute value of the largest coefficient of $f$. Hence
\[\begin{aligned}
\sum_{n=1}^N c^{\omega(f(n)/\gamma_N(n))} &\leq
\sum_{1\leq n\leq N} \mathop{\mathop{\sum_{\text{$d$ sq.-free}}}_{d\leq C \sqrt{N}}}_{d|f(n)/\gamma(n)}
c^{2 r \omega(d)} \ll_f\; N \cdot \mathop{\sum_{d\leq C \sqrt{N}}}_{p|d 
\Rightarrow p> N^{\delta(N)}} \frac{(c^{2 r} r)^{\omega(d)}}{d} \\
&\ll_{r,c} N \cdot \left(\frac{\log C \sqrt{N}}{\log N^{\delta(N)}}\right)^{
r c^{2 r}} \ll_{c,r,C} N (\log N)^{\epsilon} .\end{aligned}\]
If $\omega(f(n)/\gamma_N(n))\geq \epsilon \log \log N$, then 
$c^{\omega(f(n)/\gamma_N(n))} \geq (\log N)^{\epsilon \log c}$. We set
$c = \lceil e^{\frac{A}{\epsilon} + 1}\rceil$ and
conclude that $\omega(f(n)/\gamma_N(n))\geq \epsilon \log \log N$ for
only $O_{f,A,\epsilon}(N (\log N)^{-A})$ integers $n$ from $1$ to $N$.
\end{proof}

Now we will translate Prop.\ \ref{prop:san} into number theory. It may seem
surprising that such a thing is possible, as 
Prop.\ \ref{prop:san} assumes that the random variables it is given
are jointly independent. We will be working with the random variables
$X_p$, where $X_p = 1$ if $p$ divides a random positive integer $n\leq N$,
and $X_p=0$ otherwise; the indices $p$ range across all primes
$p\leq z$, where $z$ is such that $\log \log z\geq (1 - \epsilon) \log \log N$. 
While the variables $X_p$ are very nearly
pairwise independent, they are far from being jointly independent.
(Even if $z$ were as low as $(\log N)^2$, they would not be.)

Fortunately, the events $X_p=1$ are so rare 
($\mathbb{P}(X_p = 1) = \frac{1}{p}$) that, for a typical
$n\leq N$, the product
$d$ of all $p\leq z$ such that $X_p=1$ is at most $N^{\epsilon}$. Since
$d\leq N^{\epsilon}$, the variables $X_p$, $p|d$, are jointly independent
(up to a negligible error term). One cannot rush to conclusions, of course,
since $d$ depends on the values taken by the variables $X_p$. Nevertheless,
a careful analysis gives us the same final result as if all variables
$X_p$, $p\leq z$, were jointly independent. This procedure is not new;
it goes back in essence to Erd\H{o}s and Kac (\cite{EK}).
\begin{prop}\label{prop:agur}
Let $f\in \mathbb{Z}\lbrack x\rbrack$ be a non-constant polynomial 
irreducible over $\mathbb{Q}$.
Let the rational primes be partitioned into $\{P_j\}_{j\in J}$,
$J$ finite, so that, for every $j\in J$, we have the asymptotic
$\sum_{p\in P_j, p\leq N} 1/p \sim r_j \log \log N$,
where $\vec{r} \in (\mathbb{R}_0^+)^d$. Assume furthermore that,
for all $p\in P_j$, the equation
$f(x) \equiv 0 \mo p$ has exactly $s_j$  solutions in $\mathbb{Z}/p
\mathbb{Z}$, where
$\vec{s} \in (\mathbb{Z}_0^+)^J$. 
Let $\omega_j(n)$ be the number of divisors of $n$ in $P_j$.

Let $c_j = r_j s_j$ for $j\in J$. 
For every $\vec{x} \in (\mathbb{R}_0^+)^J$, let
\[S_{\vec{c},\vec{x}}(N) =
\{1\leq n\leq N:
(\omega_{j}(f(n)) - x_j \log \log N)\cdot
(x_j - c_j) > 0\;\;\;\; \forall j\in J\}.\]
Then, for all $\vec{x}\in (\mathbb{R}_0^+)^J$,
\[\lim_{N\to \infty} \frac{1}{\log \log N}
\log \left(\frac{1}{N} |S_{\vec{c},\vec{x}}(N)|\right)
= - I_{\vec{c}}(\vec{x}),\]
where $I_{\vec{c}}(\vec{x})$ is as in (\ref{eq:jacbo}).
\end{prop}
\begin{proof}(Cf.\ \cite{EK}, \S 4.) Let $P(z) = \prod_{p\leq z} p$.
For $d|P(z)$, let $S_{d,z}(N) = \{1\leq n\leq N: \gcd(f(n), P(z)) = d\}$.
Applying Lemma \ref{lem:spex} with $f(n)=n$, we obtain, for
$A$ arbitrarily large and $\epsilon>0$ arbitrarily small,
\begin{equation}\label{eq:fast}
\mathop{\sum_{d|P(z)}}_{d> N^{\varsigma}} |S_{d,z}(N)| = 
O_{f,A,\epsilon}\left(N (\log N)^{-A}\right),\end{equation}
where we let $z = N^{\delta(N)}$ and set $\varsigma\in (0,1)$ arbitrarily (say
$\varsigma = 1/2$). We will set and use $\epsilon$ later; for now,
it is hidden in the properties that the statement of Lemma \ref{lem:spex}
ensures for the function $\delta$ it has just defined.
 By the fundamental lemma of sieve theory
(vd., e.g., \cite{IK}, Lemma 6.3, or 
\cite{Grbook}, \S 3.3, Cor.\ 1.1) and the fact that Lemma \ref{lem:spex} 
gives us $\delta(x) = o(1/\log \log x)$,
 we have, for all
$d<N^{\varsigma}$,
\begin{equation}\label{eq:erla}
|S_{d,z}(N)| = \left(1 + O\left((\log N)^{-A}\right)\right)\cdot
\frac{N}{d} \prod_{j\in J} \mathop{\prod_{p\leq z:\, p\nmid d}}_{p\in
P_j} (1 - s_j/p) .
\end{equation}
(We use the fact
that $|\{1\leq n\leq N: \gcd(f(n),P(z))=d\}|$ equals $|\{1\leq n\leq N/d :
\gcd(f(n),P(z)/d) = 1\}|$, and estimate the latter quantity by a sieve such 
as Brun's or Rosser-Iwaniec's; we know that the sieve gives us asymptotics
with a good error term (namely,
$(\log N)^{-A}$) thanks to the fundamental lemma.)

Define the jointly independent random variables $\{X_p\}_{\text{$p$ prime}}$
as in (\ref{eq:carne}). For $d|P(z)$, let $s_{d,z}$ be the probability that
$X_p \ne 0$ for all $p|d$ and
$X_p = 0$ for all $p|P(z)/d$. By inclusion-exclusion,
$s_{d,z} = \frac{1}{d} \prod_{j\in J} \prod_{p\leq z, p\nmid d: p\in P_j} 
(1 - s_j/p)$.
Thus 
\begin{equation}\label{eq:nurmir}|S_{d,z}| = N (1 + O((\log N)^{-A})) 
\cdot s_{d,z}
\end{equation} for $d < N^{\varsigma}$, and
\begin{equation}\label{eq:sthelena}\begin{aligned}
\mathop{\sum_{d|P(z)}}_{d> N^{\varsigma}} s_{d,z} &= 
1 - \mathop{\sum_{d|P(z)}}_{d\leq N^{\varsigma}} s_{d,z} \\ &=
1 - (1 + O((\log N)^{-A})) \cdot \frac{1}{N} 
\sum_{d|P(z): d\leq N^{\varsigma}}
 |S_{d,z}|\\ &= O\left((\log N)^{-A}\right) + \frac{1}{N}
\sum_{d|P(z): d > N^{\varsigma}} |S_{d,z}| = 
O\left((\log N)^{-A}\right) ,\end{aligned}
\end{equation}
where we are using (\ref{eq:fast}) in the last line.

Since the variables $X_p$ are jointly
independent, we may apply Prop.\ \ref{prop:san}, and obtain
\[\lim_{z\to \infty} \frac{1}{\log \log z}
\log \left(\sum_{d|P(z)} \Delta(d,z)\, s_{d,z}\right) = - I_{\vec{c}}(\vec{x}) ,\]
where $\Delta(d,z) = 1$ if $
\left\{\frac{\omega_j(d)}{\log \log z}\right\}_{j\in J} \in
B_{\vec{c},\vec{x}}$ and $\Delta(d,z)=0$ otherwise.
By (\ref{eq:erla}), (\ref{eq:nurmir}) and (\ref{eq:sthelena}),
it follows that
\begin{equation}\label{eq:saus}
\lim_{N\to \infty} \frac{1}{\log \log z} \log \left(\frac{1}{N} 
\sum_{d|P(z)} \Delta(d,z)\, S_{d,z}(N)\right) = - I_{\vec{c}}(\vec{x}) 
\end{equation}
for $z = N^{\delta(N)}$ and $A$ sufficiently large, provided that
$I_{\vec{c}}(\vec{x})$ be finite. If $I_{\vec{c}}(\vec{x}) = \infty$,
we obtain (\ref{eq:saus}) with
$\lim$ replaced by $\limsup$
and $= I_{\vec{c}}(\vec{x})$ replaced by $\leq -A$. 

 Lemma
\ref{lem:spex} states that
 $|\log \delta(N)| < \epsilon \log \log N$. Thus $\log \log z > (1 - \epsilon) \log \log N$. By Lemma
\ref{lem:spex}(\ref{it:egg}), \begin{equation}\label{eq:graus}\sum_{j\in J}
 |w_j(\gcd(f(n),P(z))) - w_j(f(n))| <
\epsilon \log \log N\end{equation}
 for all but $O(N (\log N)^{-A})$ integers $n$
between $1$ and $N$.

We conclude from (\ref{eq:saus}) and (\ref{eq:graus})
that, if $I_{\vec{c}}(\vec{x})$ is finite,
\begin{equation}\label{eq:harmle}
\frac{1}{\log \log N} \log \left(\frac{1}{N} 
\sum_{n\leq N} \Delta(f(n),N)\right) = - I_{\vec{c}}(\vec{x}) +
 O_{\vec{c},\vec{x}}(\epsilon) + o_f(1),\end{equation}
where we use the fact that $I_{\vec{c}}(\vec{x})$ is continuous with
respect to the coordinate $x_j$ of $\vec{x}$ when $x_j\ne c_j$, and the fact
that the projection of $B_{\vec{c},\vec{x}}$ onto the $j$th axis is 
$\mathbb{R}$ when $x_j=c_j$. 
We let $\epsilon\to 0$ and are done.

Suppose now that $I_{\vec{c}}(\vec{x}) = \infty$. We then have 
(\ref{eq:harmle}) with $\leq -A + O_{\vec{c},\vec{x}}(\epsilon) + o_f(1)$
instead of $-I_{\vec{c}}(\vec{x}) + O_{\vec{c},\vec{x}}(\epsilon) + o_f(1)$.
We let $A\to \infty$ and $\epsilon\to 0$, and are done.
\end{proof}
It is easy to generalise Prop.\ \ref{prop:agur} so as to let
the argument $n$ of $f(n)$
range over tight sets other than the integers. (See Def.\ \ref{defn:oglo}
for the definition of a {\em tight} set.)
The means of the generalisation
will be based on a view of sieves that may be unfamiliar to some readers
and thus merits an introduction. We will use an upper-bound sieve to provide
a majorisation of the characteristic function of a tight set (such as the
primes). We will then use this majorisation as a model for the tight set,
instead of using it directly to obtain upper bounds on the number of
elements in the tight set. This model will have the virtue of being very
evenly distributed across arithmetic progressions.

We recall that an upper-bound sieve\footnote{
Take, for
example, Selberg's sieve $\lambda_d$. We are using the notation
in \cite{IK}, \S 6, and so, by $\lambda_d$, we mean the sieve coefficients, and
not the parameters (call them $\rho_d$, as in \cite{IK}) such that
$\sum_{d|n} \lambda_d = \left(\sum_{d|n} \rho_d\right)^2$. In
\cite{Grbook} and some of the older literature, the symbols $\lambda_d$
stand for what we have just denoted by $\rho_d$.
} of level $D$ is a sequence
$\{\lambda_d\}_{1\leq d\leq D}$ with $\lambda_d=1$ and
$\sum_{d|n} \lambda_d \geq 0$. Since $\lambda_d$ has support on
$\{1,2,\dotsc,D\}$, we have $\sum_{d|p} \lambda_d = 1$ for every
prime $p>D$. Thus $g(n) = \sum_{d|n} \lambda_d$ majorises the characteristic
function of $\{\text{$p$ prime}: p> D\}$. In general, if $\lambda_d$ is 
supported on $\{1\leq d\leq D: p|d \Rightarrow p\in P\}$, where $P$ is
some set of primes, $g(n) = \sum_{d|n} \lambda_d$ majorises the characteristic
function of $\{n\in \mathbb{Z}^+ : p|n \Rightarrow (p\notin P\; \vee\; p>D)\}$.

If $S\subset \mathbb{Z}^+$ is a $(P,\theta)$-tight set (vd.\ Def.\ 
\ref{defn:oglo}), then $S \cap \lbrack N^{1/2}, N\rbrack$ is contained
in $\{n\in \mathbb{Z}^+ : p|n \Rightarrow (p\notin P\; \vee\; p>N^{\delta/2})\}$,
where $\delta>0$ is as in Def.\ \ref{defn:oglo}. We set
$D = N^{\delta/2}$, and obtain that $g(n)$ majorises
the characteristic function of $S\cap \lbrack N^{1/2}, N\rbrack$. Any good
upper-bound sieve (such as Selberg's or Rosser-Iwaniec's)
 amounts to a choice of $\lambda_d$
such that $\sum_{n\leq N} g(n) \ll N/(\log N)^{\theta}$, where 
$\theta$ is the dimension of the sieving set $P$ (see 
Def.\ \ref{defn:sievdim}). Now, by Def.\ \ref{defn:oglo}, the fact that $S$ is
tight implies that 
$|S \cap \lbrack N^{1/2}, N\rbrack| \gg N/(\log N)^{\theta}$.
Thus
\begin{equation}\label{eq:asinq}
|S \cap \lbrack N^{1/2},N\rbrack| \leq \sum_{n\leq N} g(n) \ll
|S\cap \lbrack N^{1/2}, N\rbrack| .
\end{equation}
In other words, $g(n)$ is not just any majorisation of the characteristic
function of $S\cap \lbrack N^{1/2}, N\rbrack$, but a tight one, up to a 
constant factor.
\begin{prop}\label{prop:otroy}
Let $f$, $P_j$, $J$,
$\vec{r}$, $\vec{s}$, $\omega_j$ and $\vec{c}$ be
as in Prop.\ \ref{prop:agur}. Let $S\subset \mathbb{Z}$ be a 
$(P,\theta)$-tight set.
For
every $\vec{x}\in (\mathbb{R}_0^+)^J$, define
$S^*_{\vec{c},\vec{x}}(N)$ to be
\begin{equation}\label{eq:fafner}\{n\in S: 1\leq n\leq N,\,
(\omega_{j}(f(n)) - x_j \log \log N)\cdot
(x_j - c_j) > 0\;\;\, \forall j\in J\}.\end{equation}
Then, for all $\vec{x}\in (\mathbb{R}_0^+)^J$,
\begin{equation}\label{eq:omio}\limsup_{N\to \infty} \frac{1}{\log \log N}
\log \left(\frac{1}{|\{n\in S: 1\leq n\leq N\}|}
 |S^*_{\vec{c},\vec{x}}(N)|\right)
\leq - I_{\vec{c}}(\vec{x}),\end{equation}
where $I_{\vec{c}}(\vec{x})$ is as in (\ref{eq:jacbo})
 and $\delta_{\vec{x}}$ is as in Prop.\ \ref{prop:san}.
\end{prop}
The lower bound on the rate of convergence of (\ref{eq:omio})
that can be made explicit from the proof below 
depends on the constants in Def.\ \ref{defn:oglo}
(that is, on $\delta$ and the implied constant in the said definition)
but not otherwise on $(P,\theta)$. (By {\em a lower bound on the rate
of convergence} we mean a map $\epsilon\mapsto N_{\epsilon}$ 
such that the left side of (\ref{eq:omio}) is within $\epsilon$ of
the right side for all $N>N_{\epsilon}$.)

The proof of Prop.\ \ref{prop:otroy} is essentially the same as that
of Prop.\ \ref{prop:agur}; we limit ourselves to sketching the
argument again and detailing the changes.
\begin{proof}[Proof of Prop.\ \ref{prop:otroy}]
Choose an upper-bound sieve
 $\lambda_d$ of level $N^{\sigma}$,
$0 < \sigma < 1$, with
$\{p\in P : p<N^{\sigma'}\}$, $0 <\sigma' < \delta/2$, 
as its sieving set, where $P$ and $\delta$
are as in the definition of $(P,\theta)$-tight sets. 
(For example, choose $\lambda_d$ to be Selberg's sieve. See, e.g.,
\cite{IK}, \S 6.)
 The proof of Prop.\ \ref{prop:agur} goes through as before 
if one assigns the multiplicities
$\sum_{d|n} \lambda_d$ to the elements $n$ of
$\{1,2,\dotsc,N\}$, $S_{\vec{c},\vec{x}}$ and $S_{d,z}$. 
(Choose $\varsigma < 1 - \sigma$. Redo Lem.\ \ref{lem:spex} taking
into account the new multiplicities. The crucial fact is that
the natural estimates for $\sum_{1\leq n\leq N: r|n} \sum_{d|n} \lambda_d$
($r$ given)
have very good error terms. The irreducibility of $f$ helps us in so far
as $f(x)\equiv 0 \mo p$, $p|x$ are both true for a finite number of
primes $p$, if for any.) We obtain the statement of Prop.\ \ref{prop:agur},
with
\begin{equation}\label{eq:hug}\lim_{N\to \infty} \frac{1}{\log \log N} \log 
\left(\frac{1}{\sum_{1\leq n\leq N} \sum_{d|n} \lambda_d} 
 |S_{\vec{c},\vec{x}}|\right) = -I_{\vec{c}}(\vec{x})\end{equation}
as the result, and $S_{\vec{c},\vec{x}}$ counting $n$ with the multiplicity
$\sum_{d|n} \lambda_d$.
Since $\sum_{d|n} \lambda_d$ majorises the characteristic function of
$S \cap \lbrack N^{1/2}, N\rbrack$, we have 
$|S_{\vec{c},\vec{x}}^*| \leq |S_{\vec{c},\vec{x}}| + N^{1/2}$. At the same time, as in (\ref{eq:asinq}),
$\sum_{n\leq N} \sum_{d|n} \lambda_d \ll |S \cap \lbrack N^{1/2}, N\rbrack|$.
Hence (\ref{eq:hug}) implies (\ref{eq:omio}).
\end{proof}
For $S$ equal to the set of all primes, we could replace
$\limsup$ and $\leq$ in (\ref{eq:fafner}) by $\sup$ and $=$ through
an appeal to Bombieri-Vinogradov. However, we shall not need such an 
improvement.
\section{Proof of the main theorem and immediate consequences}\label{sec:pmt}
Using the results in \S \ref{sec:largdev}, we will now show that,
if $t p^{k} = f(n)$ for some $n\in S\cap \lbrack 1,N\rbrack$, some
prime $p$ and some $t$ not much larger than $N$, then either $n$ is
atypical or $t$ is atypical. Since ``atypical'' means ``rare'', we 
conclude, counting either $n$'s or $t$'s, that few $n\in S\cap \lbrack 1,N\rbrack$ satisfy $t p^{k} = f(n)$ for some $p$ prime and some integer
$t$ not much larger than $N$.

For this argument to yield anything of use to us, we must make it quantitative
and rather precise. It is here that entropies come into play, as they
appear in the exponents of expressions for the probabilities of unlikely
events.
\begin{prop}\label{prop:knut}
Let $f\in \mathbb{Z}\lbrack x\rbrack$ be a polynomial of degree $d\geq 3$
irreducible over $\mathbb{Q}$. Let $S$ be a 
$(P,\theta)$-tight set of integers. Define $J$ to be the set
of conjugacy classes $\langle g\rangle$ of $\Gal_f$. Let
$\vec{c} = \left\{\frac{|\langle g\rangle|}{|\Gal_f|}\right\}_{\langle g\rangle
\in J}$ and $\vec{c}' = \left\{\frac{|\langle g\rangle|}{|\Gal_f|} 
\lambda_{\langle g\rangle} \right\}_{\langle g\rangle \in J}$,
where $\lambda_{\langle g\rangle}$ is the number of roots of
$f(x)=0$ fixed by $g\in \Gal_f$. 
 Define
\begin{equation}\label{eq:odi}
\gamma = \min_{\vec{x}\in X} \left(\max\left(\frac{d-1}{d} + \frac{1}{d}
I_{\vec{c}}(\vec{x}),\; \theta + I_{\vec{c}'}(\vec{x})\right)\right) ,
\end{equation}
where $X = \vartimes_{j\in J} \lbrack \min(c_j,c_j'), \max(c_j,c_j')\rbrack$
and $I_{\vec{c}}(\vec{x})$, $I_{\vec{c}'}(\vec{x})$ are as in (\ref{eq:jacbo}).
Then
\begin{equation}\label{eq:schuby}
\frac{1}{N} |\{n\in S\cap \lbrack 1,N\rbrack : \text{$\exists p$ prime,
$p\geq N^{\epsilon}$ such that $p^{d-1}|f(n)$}\}|\end{equation}
is $O((\log N)^{-\gamma + \epsilon})$ 
for every $\epsilon>0$. The implied constant depends only on
$f$, $\theta$, $\epsilon$, and the constants in Def.\ \ref{defn:oglo} for
the given set $S$.
\end{prop}
\begin{proof}
There is a $t>0$ such that $X\subset \lbrack 0,t\rbrack^J$ and 
$I_{\vec{c}'}(\vec{x}) > \gamma$ for all
elements $\vec{x}$ of $\{\vec{x} \in (\mathbb{R}_0^+)^J : x_j \geq t\;
\text{ for some $j\in J$}\}$. Let $Y = \lbrack 0, t\rbrack^J$. Cover
$Y$ by all sets of the form $B_{\vec{c},\vec{x}}$ (see (\ref{eq:tuastort}))
with $\vec{x}$ such that $\frac{d-1}{d} + \frac{1}{d} 
I_{\vec{c}}(\vec{x})\geq \gamma - \epsilon$, and all sets of
the form $B_{\vec{c}',\vec{x}}$ with $\vec{x}$ such that
$\theta + I_{\vec{c}'}(\vec{x})\geq \gamma - \epsilon$. 
Such sets
form a cover of $Y$ by (\ref{eq:odi}). 
Since $Y$ is compact
and all sets 
$B_{\vec{c},\vec{x}}$, $B_{\vec{c}',\vec{x}}$ in the cover are open, there
is a finite subcover $\mathscr{B}$; we may choose one such finite subcover
in a way that depends only on $\vec{c}$, $\vec{c}'$, $d$ and $\theta$,
and thus only on $f$ and $\epsilon$. 
Write $\mathscr{B} = 
\bigcup_{x\in X} 
B_{\vec{c},\vec{x}} \; \cup \; \bigcup_{x\in X'} B_{\vec{c}',\vec{x}}$.

Define $\mathscr{B}'$ to be the union of $\mathscr{B}$
and the collection of all sets $B_{\vec{c}',\vec{x}}$ with $\vec{x}$
such that $x_k = t$ for some $k\in J$ and $x_j = c_j'$ for
all $j\ne k$. Then $\mathscr{B}'$ is a cover of $(\mathbb{R}_0^+)^J$.
In particular, for every $n\in S\cap \lbrack 1, N\rbrack$ such that
$p^{d-1}|f(n)$ for some prime $p\geq N^{\epsilon}$, we have 
$\{w_j(f(n))/\log \log N\}_j\in B$ for some 
$B$ in $\mathscr{B}'$. 
The set $B$ may be of type 
$B = B_{\vec{c},\vec{x}}$ or
$B = B_{\vec{c}',\vec{x}}$.
In the latter case,
(\ref{eq:fafner}) holds with $\vec{c}'$ instead of $\vec{c}$,
and so, by Prop.\ \ref{prop:otroy}, $n$ belongs to a set
$S_{\vec{c}',\vec{x}}^*(N)$ whose cardinality is bounded above by a 
constant times
\begin{equation}\label{eq:ontime}
(\log N)^{-I_{\vec{c}'}(\vec{x}) + \epsilon} \cdot
|S \cap \lbrack 1,N\rbrack| \ll_{\delta} 
N (\log N)^{-I_{\vec{c}'}(\vec{x}) -
\theta + \epsilon} \ll N (\log N)^{-\gamma + 2 \epsilon} ,
\end{equation}
where $\delta$ is as in Def.\ \ref{defn:oglo}. (Here we are assuming,
as we may, that $N$ is larger than some constant depending only on
$f$, $\epsilon$, $\delta$, and the implicit
constant in Def.\ \ref{defn:oglo}.
As in further applications of Prop.\ \ref{prop:agur} and Prop.\
\ref{prop:otroy}, we define the sets $P_{\langle g\rangle}$
to consist of the primes $p\nmid \Disc(f)$ with specified 
Frobenius element $\Frob_p = \langle g\rangle \in J$; we 
put each prime $p|\Disc(f)$
in its own exceptional set $P_{0,p}$. The
exceptional sets will have no influence
on the bounds. The densities $\lim_{N\to \infty}
\frac{1}{\log \log N} \sum_{p\in P_j, p\leq N} 1/p$
of the sets $P_j$ are given by the Chebotarev density theorem.)

Consider the other possibility, namely, that $\{w_j(f(n))\}$ is in a set
$B_{\vec{c},\vec{x}}$. Let $p$ be a prime $\geq N^{4/5}$ such that
$p^{d-1} | f(n)$, and define $r = f(n)/p^{d-1}$. Suppose first that
$|r| > N (\log N)^{-\alpha}$, where $\alpha$ will be set later.
Then $p^{d-1} \ll N^d/|r| < N^{d-1} (\log N)^{\alpha}$, and so
$p\ll N (\log N)^{\alpha/(d-1)}$. For every $p$, the
number of positive integers $n\leq N$ with $p^{d-1} | f(n)$ is
$\leq O_f(\lceil N/p^{d-1} \rceil)$. Since $d\geq 3$,
\begin{equation}\label{eq:desh}\begin{aligned}
\sum_{N^{\epsilon} \leq p \leq N (\log N)^{\alpha/(d-1)}} \left\lceil 
\frac{N}{p^{d-1}} \right\rceil&\leq
\sum_{N^{\epsilon} \leq p \leq N (\log N)^{\alpha/(d-1)}} \left(
\frac{N}{p^{d-1}} + 1\right)\\
&\ll N^{1 - \epsilon} + N ( \log N)^{\alpha/(d-1) - 1} .
\end{aligned}\end{equation}
Choose $\alpha = \frac{d-1}{d} (1 - I_{\vec{c}}(\vec{x}))$. Then
$\frac{\alpha}{d-1} - 1 = - \left(\frac{d-1}{d} + \frac{1}{d} I_{\vec{c}}(
\vec{x})\right) < - \gamma + \epsilon$, and thus the contribution of all $n$ with
$r p^{d-1} = f(n)$, $|r| > N (\log N)^{-\alpha}$, is at most 
$O_f(N (\log N)^{-\gamma +\epsilon})$. 

Now take the remaining possibility, namely, 
$|r| \leq N (\log N)^{-\alpha}$. Choose an integer divisor
$k>1$ of $d-1$. Let $w$ be the product of $q^{\lfloor v_q(r)/k\rfloor\cdot k}$
over all primes $q$ dividing $k \Disc(f)$. Let $y = w p^{(d-1)/k}$,
$t = r/w^k$. Then $t y^k = f(n)$. Moreover, $|t| \leq N (\log N)^{-\alpha}$,
$\gcd(t,(k \Disc(f))^{\infty}) \leq (k \Disc(f))^k = O_{f,k}(1)$,
and the number of prime divisors $\omega(y)$ of $y$ is also $O_{f,k}(1)$.

We may assume without loss of generality that the leading coefficient of $f$
is positive, and thus $f(n)$ will be positive for $n$ larger than some
constant $O_f(1)$. Since $y$ is also positive, $t$ is positive as well.
For every $j\in J$, we know that $\omega_j(y) = O_f(1)$ and 
$\omega_j(f(n)) - \omega_j(y^k) \leq \omega_j(f(n)/y^k)
= \omega_j(t)\leq \omega_j(f(n))$. Hence 
\[\omega_j(f(n)) - O_f(1) \leq \omega_j(t) \leq \omega_j(f(n)) .\]
Define $\vec{x}' \in (\mathbb{R}_0^+)^J$ by $x_j' = c_j + (1 - \epsilon)
(x_j - c_j)$.
Then $I_{\vec{c}}(\vec{x}') \geq I_{\vec{c}}(\vec{x}) - O_f(\epsilon)$
and $\{w_j(t)\}_j \in B_{\vec{c}}(\vec{x}')$ for $n$ larger than some
constant $O_f(1)$. We may ignore all $n$ smaller than $O_f(1)$, as they
will contribute at most $O_f(1)$ to the final bound on (\ref{eq:schuby}).

We apply Prop.\ \ref{prop:agur} with
$f(x) = x$ and $N (\log N)^{-\alpha}$ instead of $N$.
We obtain that 
$t$ lies in a set $S' = S_{\vec{c},\vec{x}'}(N (\log N)^{-\alpha})$ of 
cardinality at most
\begin{equation}\label{eq:cardio}
(\log N)^{-I_{\vec{c}}(\vec{x}') + \epsilon} \cdot |\mathbb{Z}
\cap \lbrack 1, N (\log N)^{-\alpha}\rbrack|
= N (\log N)^{-I_{\vec{c}}(\vec{x}) - \alpha + \epsilon},
\end{equation}
provided that, as we may assume, $N$ is larger than some constant
$O_f(1)$.

Our task is to bound, for each $t\in S' \cup \lbrack 1, N (\log N)^{-\alpha}
\rbrack$, how many solutions $(n,y) \in (\mathbb{Z}^+)^2$ with
$n\leq N$ the equation $t y^k = f(n)$ has. (We are also given that
$\gcd(t,(k \Disc(f))^{\infty})$ is bounded above by $O_{f,k}(1)$.)
Of $S'$ we need only remember that it is a subset of $\mathbb{Z}^+$
with cardinality at most (\ref{eq:cardio}). 

Let $\delta(x)$ be as in Lemma \ref{lem:spex} with $A$ equal to
$\gamma$ (or greater).
Assume that $n$ is such that
(\ref{it:oyop}) and (\ref{it:egg}) in Lemma \ref{lem:spex}
both hold. (By the said Lemma, we are thereby excluding at most
$O_{f,A,\epsilon}(N (\log N)^{-A})$ values of $n$.) 
We may also assume that $f(n)$
has at most $O_A(\log \log n)$ prime divisors and exclude
thereby at most $O(N (\log N)^{-A})$ values of $n$.
We may also assume that
$n>N^{1-\epsilon}$ (and exclude an additional set of $N^{1-\epsilon}$ 
values of $n$).
Apply Prop.\ \ref{prop:notung} with $t_0=t/t_1$,
where $t_1$ is the product of
all primes $p|f(n)$ such that $p\leq N^{\delta(N)}$.
We obtain that there are at most 
$O_{f,k,\epsilon}((\log N)^{O_{f,k,A}(\epsilon)})$ possible values
of $n$ for every value of $t$. 

Since the number of values of $t$ under consideration
is bounded by (\ref{eq:cardio}) and $\alpha$ has been chosen so that
$-I_{\vec{c}}(\vec{x}) - \alpha = 
-\frac{d-1}{d} - \frac{1}{d} I_{\vec{c}}(\vec{x}) \leq
- \gamma + O_f(\epsilon)$, we conclude that there are at most
\begin{equation}\label{eq:anont}
O(N (\log N)^{-A}) + 
 O_{f,k,\epsilon}(|S'| (\log N)^{O_{f,k,A}(\epsilon)})
\ll_{f,k,\epsilon} N (\log N)^{-\gamma + O_{f,k}(\epsilon)} \end{equation}
solutions $(n,y) \in (\mathbb{Z}^+)^2$ with
$n\leq N$ to a given equation $t y^k = f(n)$ with 
$t\in S' \cup \lbrack 1, N (\log N)^{-\alpha}$ and
$\gcd(t,(k \Disc(f))^{\infty}) = O_{f,k}(1)$.

We add the bounds
(\ref{eq:ontime}), (\ref{eq:cardio}), (\ref{eq:desh}) and (\ref{eq:anont})
over all elements of the cover $\mathscr{B}'$. Since the cardinality of
the cover depends only on $f$, $\theta$ and $\epsilon$, we are done.
\end{proof}

Proposition \ref{prop:knut} was the ultimate purpose of all of the work
that came before it. The following lemma is far softer.

\begin{lem}\label{lem:knotty}
Let $f\in \mathbb{Z}\lbrack x\rbrack$ be a polynomial of degree $d\geq 1$
irreducible over $\mathbb{Q}\lbrack x\rbrack$. Let $k\geq \max(d,2)$.
Then there is a $\delta>0$ depending only on $d$ such that,
for every $D>1$,
\begin{equation}\label{eq:romd}
|\{1\leq n\leq N: \text{$\exists p$ prime,
$p\geq D$ \text{\; s.t.\ \;} $p^{k}|f(n)$}\}|\ll_{f} N (N^{-\delta} +
D^{-(k-1)}) .\end{equation}
\end{lem}
\begin{proof}
Suppose $k>d$. If
 $p^k|f(n)$, then $p \ll_f N^{d}{k}$.
The number of positive integers $n\leq N$
such that $p^k|f(n)$ for some prime 
$D \leq p \ll_f N^{d/k}$
can be shown to be
$O_f(N^{d/k})$ by the same simple argument as in
(\ref{eq:desh}). We set $\delta = 1 - \frac{d}{k}$ and are done.

Suppose now $k=d$. Again as in (\ref{eq:desh}), the number of
integers $n\in \lbrack 1,N\rbrack$ such that $p^k | f(n)$ for
some prime $D\leq p \ll_f N^{1 - \epsilon/k}$, $\epsilon>0$, is at most
$O_f(N^{1-\frac{\epsilon}{k}} + N D^{-(k-1)})$.
If $p^k|f(n)$ for some prime $p> N^{1-\epsilon/k}$, then $r p^k = f(n)$,
where $r$ is an integer with $|r| \ll_f N^{\epsilon}$. Thus, we need only
show that, for every integer $r$ with $|r|\ll_f N^{\epsilon}$,
\begin{equation}\label{eq:oglogo}
|\{1\leq n\leq N: \exists \text{ $p$ prime such that $r p^k = f(n)$}\}|
\ll N^{1 - 2 \epsilon}\end{equation}
provided that $\epsilon$ be sufficiently small.
This is an easy bound; a much stronger
one (viz., $\ll N^{1/k}$ instead of $\ll N^{1 - 2 \epsilon}$)
 follows immediately from \cite{BP}, Thm.\ 5. Set $\delta = \epsilon/k$.
\end{proof}
\begin{remark}
We can actually replace the bound
$\ll N^{1 - 2\epsilon}$ in the right side of (\ref{eq:oglogo})
by $\ll N^{\epsilon'}$, with $\epsilon'>0$ arbitrarily small, provided
that $\deg(f)>1$. (If $\deg(f)=1$, we have instead 
$\ll (N/r)^{1/k}/\log((N/r)^{1/k})$,
which can just as easily be proven as be proven best: let $p$ vary, and define
$n$ in terms of $r$ and $p$.) We may proceed as follows:
\begin{enumerate}
\item If $\deg(f)=2$ and $k=2$, then the number of points $(x,y)\in
(\mathbb{Z}\cap \lbrack 1, N\rbrack)^2$ on $r y^2 = f(x)$ is
$O_f(\log N)$.
 This is a classical result of
Estermann's (\cite{Es}, p.\ 654 and p.\ 656). 
(Reduce the problem to the case where $f(x)$
is of the form $x^2 + l$, $l\ne 0$, by a change of variables over $\mathbb{Q}$.
Then count the solutions $(x,y)\in (\mathbb{Z} \cap \lbrack 1,N\rbrack)^2$
to $r y^2 - x^2 = l$; they are bounded by $O_l(\log N)$ because the group
of units of $\mathbb{Q}(\sqrt{r})$ is of rank $1$.)
\item If $\deg(f)>2$ or $k>2$, the genus of $C:r y^k = f(x)$ is positive.
Bound the rank of $C(\mathbb{Q})$ by Lem.\ \ref{lem:angr} (generalised
so as to remove the assumption $p\nmid d$; see the comment after
the statement of Prop.\ \ref{prop:notung}). Bound the number of integer
solutions to $r y^k = f(x)$ with 
$X^{(1 - \epsilon'') \sigma} < x\leq X^{\sigma}$,
$\sigma\leq 1$, as in Prop.\ \ref{prop:notung}; the auxiliary divisor
$t_0$ is not needed, as we do not aim at estimates as delicate as before.
We obtain a bound of $O_{f,k}((\log N)^c)$, $c>0$ fixed, for the number
of integer points with $x$ in the said range. Vary $\sigma$ as needed.
\end{enumerate}
We obtain Lemma \ref{lem:knotty} with $\delta = \frac{d}{k+1} - \epsilon$,
where $\epsilon>0$ is arbitrary. (The implied constant in (\ref{eq:romd})
then depends on $\epsilon$.)

Alternatively, we could bound the number of rational solutions to 
$r y^k = f(x)$ of height $O(\log N)$ by Cor.\ 4.3 and Lem.\ 4.4 of
\cite{Hsq} and Prop.\ 3.6 of \cite{HV}, say, and then bound the number of 
integer solutions by the number of rational solutions. The resulting
bound
would still be $O(N^{\epsilon})$ on the average of $r$, and so
we would still get $\delta = \frac{d}{k+1} + \epsilon$, $\epsilon>0$.

 As it happens, Lemma \ref{lem:knotty} in its presently stated form 
(that is, with $\delta>0$ unspecified) is all
we shall need; even in the explicit result for prime arguments
(Prop.\ \ref{prop:primos}), the error terms would not be affected
by any improvements on Lemma \ref{lem:knotty}. The argument just sketched
in this remark was well within the reach of previously known techniques;
it has been included only for completeness.
\end{remark}

We are now ready to prove the main theorem. Given
Prop.\ \ref{prop:knut}, what remains is quite straightforward.
\begin{proof}[Proof of Main Theorem]
Our main task is to show that, for every $\epsilon>0$,
\begin{equation}\label{eq:coffee}
|\{1\leq n\leq N: \text{$\exists p$ prime, $p\geq N^{\epsilon}$ such that
$p^k | f(n)$}\}| = o_{f,S,\epsilon}(N) .\end{equation}
Let $f = c f_1^{r_1} \dotsb f_l^{r_l}$, where $c\in \mathbb{Z}$,
$r_i < f_k$ and the $f_i$'s are irreducible polynomials in 
$\mathbb{Z}\lbrack x\rbrack$ coprime to each other. Then, for $q$ larger than
a constant, we may have $q^k|f(p)$ only if $q^k | f_i^{r_i}$ for some
$i\in \{1,2,\dotsc,l\}$. Note that $q^k|f_i^{r_1}$ implies 
$q^{k_i}|f_i$, where $k_i$ equals $\lceil \frac{k}{r_i}\rceil$, which,
by the assumption in the statement of the theorem, is at least $\deg(f_i)-1$.
Thus, for the purpose of proving (\ref{eq:coffee}), we may assume that
$f$ is irreducible and $k\geq \deg(f) - 1$.

If $k\geq \deg(f)$, then (\ref{eq:coffee}) follows immediately from Lemma
\ref{lem:knotty}. Suppose $k = \deg(f) - 1$. By Prop.\ 
\ref{prop:knut}, we need only check that $\gamma$ as defined in
\ref{eq:odi} is greater than $\theta$. Since
$I_{\vec{c}}(\vec{x})$ is continuous on $\vec{x}$ in the domain on which
it is finite, it is enough to check that
 $\frac{k}{k+1} + 
\frac{1}{k+1} I_{\vec{c}}(\vec{c}') > \theta$, 
as it will then follow that
$\frac{k}{k+1} + 
\frac{1}{k+1} I_{\vec{c}}(\vec{x}) > \theta + \epsilon'$ for
some $\epsilon'>0$ and any $\vec{x}$ in
some open neighbourhood of $\vec{c}'$, and, by (\ref{eq:jacbo}),
$\theta + I_{\vec{c}'}(\vec{x}) > \theta + \epsilon''$, $\epsilon''>0$,
outside that neighbourhood.

We must, then, show that $\frac{k}{k+1} + 
\frac{1}{k+1} I_{\vec{c}}(\vec{c}') > \theta$. Now,
\[\begin{aligned}
I_{\vec{c}}(\vec{c}') &= 1 - \sum_{\langle g\rangle} c_{\langle g\rangle}' +
\sum_{\langle g\rangle} c_{\langle g\rangle}' \log \frac{c_{\langle g\rangle}'}{c_{\langle g\rangle}} = \sum_{\langle g \rangle} c_{\langle g\rangle}'
\log \lambda_{\langle g\rangle}\\
&= \frac{1}{|\Gal_f|} \sum_{\langle g\rangle} |\langle g\rangle| 
\lambda_{\langle g\rangle} \log \lambda_{\langle g\rangle} =
\frac{1}{|\Gal_f|} \sum_g \lambda_g \log \lambda_g = I_f, 
\end{aligned}\]
where we use the fact that $\sum c_j' = 1$ (by the Cauchy-Frobenius Lemma,
or, as it is incorrectly called, Burnside's Lemma; see \cite{Ne}).
By one of the assumptions in the statement of the present theorem,
$I_f > (k+1) \theta - k$. Thus, $\frac{k}{k+1} +
\frac{1}{k+1} I_f > \theta$. We are done proving (\ref{eq:coffee}).

Since $S$ is $(P,\theta)$-tight, we have, for every $p\leq N^{\epsilon}$,
\[|\{1\leq n\leq N : p^k|f(n)\}| = O_f(N/p^k + 1),\]
where we use an upper-bound sieve with sieving set $P\setminus
(P\cap \{p\})$ to bound the cardinality on the left. 
(The bound on the right is attained by the definition of
$(P,\theta)$-tightness.)
Thus, for any $z>0$, \begin{equation}\label{eq:peer}
|\{1\leq n\leq N: \text{$\exists p$ prime,
$z<p\leq N^{\epsilon}$ s.t.\ $p^k|f(n)$}\}| = 
O_{f,\epsilon}(N/z^{k-1} + N^{\epsilon}) . \end{equation}
Let $A_{f,k,z}(N)$ be the set of integers $n\in S\cap \lbrack 1,N\rbrack$
such that $p^k\nmid f(n)$ for every $p\leq z$. Since $S$ is
predictable (see Def.\ \ref{def:boul}), 
\begin{equation}\label{eq:gynt}\begin{aligned}
\lim_{N\to \infty} \frac{|A_{f,k,z}(N)|}{|S\cap \lbrack 1,N\rbrack|} 
&= \mathop{\sum_{m\leq 1}}_{p|m \Rightarrow p\leq z} \mu(m)\cdot
\lim_{N\to \infty} \frac{|\{n\in S\cap \lbrack 1,N\rbrack: m^k|f(n)\}|}{
|S\cap \lbrack 1,N\rbrack|} \\
&= 
\mathop{\sum_{m\leq 1}}_{p|m \Rightarrow p\leq z} \mu(m)
\mathop{\sum_{0\leq a<m^k}}_{f(a)\equiv 0 \mo m^k} \rho(a,m^k)
 ,\end{aligned}\end{equation}
where the rate of convergence depends on $z$, which is here held fixed.
Let $A_{f,k}(N)$ be the set of integers $n\in S\cap \lbrack 1,N\rbrack$ such that
$f(n)$ is free of $k$th powers. By (\ref{eq:coffee}), (\ref{eq:peer}) and
(\ref{eq:gynt}),
\[\begin{aligned}
|A_{f,k}(N)| &= |S \cap \lbrack 1,N\rbrack|\cdot
\mathop{\sum_{m\leq 1}}_{p|m \Rightarrow p\leq z} \mu(m)
\mathop{\sum_{0\leq a<m^k}}_{f(a)\equiv 0 \mo m^k} \rho(a,m^k)\\
&+ o_{z,S}(N) + O_{f,S,\epsilon}(N/z^{k-1} + N^{\epsilon}) + o_{f,\epsilon}(N)
\end{aligned}\]
for every $z$. Choose $\epsilon = \frac{1}{5}$ (say). We let $z$ go
to infinity with $N$ as slowly as needed, and conclude that 
\begin{equation}\label{eq:koechlin}
|A_{f,k}(N)| = |S \cap \lbrack 1,N\rbrack|\cdot
\mathop{\sum_{m\leq 1}}_{p|m \Rightarrow p\leq z} \mu(m)
\mathop{\sum_{0\leq a<m^k}}_{f(a)\equiv 0 \mo m^k} \rho(a,m^k)
+ o_{f,S}(N) .\end{equation}
\end{proof}
\begin{remark}
The implied constant in (\ref{eq:koechlin}) depends only on $f$, on the 
constants in Def.\ \ref{defn:sievdim} and on the rate of convergence in
(\ref{eq:taois}) for the given set $S$ (as a function of $m$).
\end{remark}
\begin{proof}[Proof of Theorem \ref{thm:cea}]
The primes are a $(P,\theta)$-tight set with $\theta=1$ and $P$ equal
to the set of all primes; the constant $\delta$ in Def.\ \ref{defn:oglo}
is $1$. By the prime number theorem, the expression $\rho(a,m)$ in
(\ref{eq:taois}) equals $\frac{1}{\phi(m)}$ if $\gcd(a,m)=1$ and $0$
if $\gcd(a,m)\ne 1$. In particular, the primes are predictable.
Apply the main theorem. Since $\phi(m)$ is multiplicative, the expression
(\ref{eq:bormo}) equals (\ref{eq:boras}).
\end{proof}

\begin{proof}[Proof of Corollary \ref{cor:orovil}]
Let $c$ be the greatest common divisor of the
coefficients of $f$. For every $p>d+1$, the equation
$\frac{1}{p^{v_p(c)}} f(x) \equiv 0 \mo p$ has at most
$d < p-1$ solutions in $\mathbb{Z}/p \mathbb{Z}$. Hence
$\frac{1}{p^{v_p(c)}} f(n) \not\equiv 0 \mo p$ for some integer
$n$ not divisible by $p$. 
Clearly $f(n)\not\equiv 0 \mo p^{v_p(c)+1}$. Since
$v_p(c) < k$, we conclude that 
$f(x)\not\equiv 0 \mo p^{k}$ has a solution in
$(\mathbb{Z}/p^k \mathbb{Z})^*$ for every $p>d+1$. 
We are given, by the assumption in the statement, that $f(x)\not\equiv 0
\mo p^k$ has a solution in $(\mathbb{Z}/p^k \mathbb{Z})^*$ for every
$p\leq d+1$ as well. We obtain that no factor of
(\ref{eq:boras}) is $0$.

For $p$ sufficiently large, $f(x)\equiv 0 \mo p^k$ has at most
$k$ solutions in $(\mathbb{Z}/p^k \mathbb{Z})^*$, by Hensel's Lemma;
thus
\[1 - \frac{\rho_{f,*}(p^k)}{p^k - p^{k-1}} \geq
1 - \frac{k}{p^{k} - p^{k-1}} \geq 1 - \frac{2 k}{p^2},\]
and so we see that (\ref{eq:boras}) does not converge to $0$.
Apply Thm.\ \ref{thm:cea}.
\end{proof}
Now let us show that,
as was remarked at the end of \S \ref{subs:prima}, the entropy
 $I_f$ is greater than $1$ for all normal
polynomials $f$ of degree $\geq 3$, and, 
in particular, for all $f$ with $\Gal_f$ abelian and $\deg(f)\geq 3$.

Let $f\in \mathbb{Z}\lbrack x\rbrack$ 
be a polynomial of degree $d\geq 3$. Suppose $f$ is normal, i.e.,
 $|\Gal_f| = d$. By the Cauchy-Frobenius Lemma, 
$\frac{1}{|\Gal_f|} \sum_{g\in \Gal_f} \lambda_g = 1$. At the same time,
$\lambda_{e} = d$ for the identity element $e\in \Gal_f$. Hence
$\lambda_{g} = 0$ for every $g\in \Gal_f$ other than the identity.
So,
\begin{equation}\label{eq:squiggle}
I_f = 
\frac{1}{|\Gal_f|} \mathop{\sum_{g\in \Gal_f}}_{\lambda_g\ne 0}
\lambda_g \log \lambda_g = \frac{1}{d} \cdot d \log d = \log d > 1.
\end{equation}
Note that $\log d$ is the largest entropy a polynomial
 of degree $d$ can have.

A transitive abelian group on $n$ elements has order $n$ (see,
e.g., \cite{Sc}, 10.3.3--10.3.4). Thus, every polynomial $f$
with $\Gal_f$ abelian is normal, and, by the above, its entropy $I_f$
is $\log(\deg(f))>1$.

Lastly, let us compute the entropy of $f$ with $\Gal_f = S_n$ and
$n$ large. (A generic polynomial of degree $n$ has Galois group 
$\Gal_f = S_n$.) Let us define the random variable $Y_n$ to be
the number of fixed points of a random permutation of 
$\{1,2,\dotsc,n\}$. Then, for $f$ irreducible with $\Gal_f = S_n$,
\[I_f = \sum_{k=1}^n \Prob(Y_n = k) k \log k.\] It is easy to show that
the distribution of $Y_n$ tends to a Poisson distribution as $n\to \infty$;
in fact, by, say, \cite{AT}, pp.\ 1567--1568,
\[\max_{0\leq k\leq n} \left| \Prob(Y_n = k) - \frac{e^{-1}}{k!} \right|
\leq \frac{2^{n+1}}{(n+1)!} .\] Thus,  
$I_f = \sum_{k=1}^n \frac{e^{-1}}{k!} k \log k + o(1)$, and so
\begin{equation}\label{eq:onyo}
\lim_{n\to \infty} I_f 
= \sum_{k=1}^{\infty}
\frac{e^{-1}}{k!} k \log k .
\end{equation}
Numerically, $\sum_{k=1}^{\infty} \frac{e^{-1}}{k!} k \log k =
0.5734028\dotsc$. Since this is less than 1, the conditions of
Thm.\ \ref{thm:cea} are not fulfilled for $f$ with $\deg(f)=n$,
$\Gal_f = S_n$, $n$ large; some simple numerics suffice to show
the same (namely, $I_f<1$) for $f$ with $\deg(f)=n$, $\Gal_f = S_n$,
$n$ small. This is unfortunate, as a generic polynomial of degree $n$
has Galois group $S_n$.
\section{Rates of convergence and error terms}\label{sec:rater}
We now wish to bound the error terms implicit in the various cases of
the main theorem discussed in \S \ref{sec:introdu}. We must first compute
 the quantity $\gamma$ defined in (\ref{eq:odi}).

This computation may seem familiar to those who have seen the theory of
large deviations being used in hypothesis testing. Let us go through
a simple example of such a use. Some believe that the variable $X$
follows a certain distribution, centred at $a$, say; others believe
it follows another distribution, centred at $b$. Both parties agree to
fix a threshold $c$ (with $a<c<b$) and take $n$ observations of the variable $X$.
If the sample mean $S_n = \frac{1}{n} (X_1 + \dotsc + X_n)$ turns out to be 
less than $c$, the contest will have been decided in favour of the
distribution centred at $a$; if $S_n>c$, the distribution centred at $b$
will be held to be the correct one. The question is: where is the
best place to set the threshold? That is, what should $c$ be?

Denote by $\Prob_a(E)$ the probability of an event $E$ under the assumption
that the distribution centred at $a$ is the correct one, and by
$\Prob_b(E)$ the probability if the distribution centred at $b$ is
the correct one. Then we should set $c$ so that
$\max(\Prob_a(S_n>c),\Prob_b(S_n<c))$ 
is minimal; that way, the likelihood
of resolving the contest wrongly will be minimised. (We are making no
{\em a priori} assumption as to the likelihood of either party being correct.)
This minimum will usually be
attained when $\Prob_a(S_n>c) = \Prob_b(S_n<c)$.
Actually computing $c$ is a cumbersome task; it is rare that
there is a closed expression either for the minimum of
$\max(\Prob_a(S_n>c),\Prob_b(S_n<c))$ 
 or for the $c$ for which it is attained.

In our context, we have that any value of $d$ in $d y^2 = f(p)$ will
be unlikely either as an integer or as a value of $f(p)$ divided
by the square of some prime. We must set a threshold of some sort and be able
to say that, if $d$ falls under it (in some sense), it must be unlikely
as an integer, and, if it goes over it, it must be unlikely as a value
of $f(p)$ divided by the square of a prime. It will be best to set the
threshold so that the maximum of the two likelihoods will be minimised.
(Matters are complicated by the facts that, in our problem,
 one of the distributions starts ``$(\log N)^{-\theta}$ ahead''
($\theta=1$ in the case of prime argument); thus we have (\ref{eq:odi})
 instead of $\min_c \max(\Prob_a(S_n>c),\Prob_b(S_n<c))$.)
The minimum of the maximum of the two likelihoods can usually be attained
only when the two likelihoods are equal; we have to minimise
them on the surface on which they are equal.
(We will be working in several dimensions. Thus,
the fact that the two likelihoods are equal defines a surface.)
  
Again, it is difficult to give closed expressions for all constants, but
we shall always be able to compute all exponents -- and, in some particular
cases, fairly simple expressions can in fact be found; see the third
note after the proof of Proposition \ref{prop:primos}.

\begin{lem}\label{lem:mamsol}
Let $\vec{c}, \vec{c}' \in \left(\mathbb{R}_0^+\right)^J$,
$J$ finite. Define
\begin{equation}\label{eq:empan}
g(\vec{x}) = \max\left(\frac{d-1}{d} + \frac{1}{d}
I_{\vec{c}}(\vec{x}),\; \theta + I_{\vec{c}'}(\vec{x})\right),\end{equation}
where
$I_{\vec{c}}(\vec{x})$, $I_{\vec{c}'}(\vec{x})$ are as in (\ref{eq:jacbo}).
Define 
\[X = \vartimes_{j\in J} \lbrack \min(c_j,c_j'), \max(c_j,c_j')\rbrack .\]
Then the minimum of $g(\vec{x})$ on $X$ is attained when and only when
\begin{equation}\label{eq:horna}
x_j = \begin{cases}
 c_j^{\alpha} (c_j')^{1-\alpha} &\text{if $c_j, c_j'\ne 0$,}\\
0 &\text{if $c_j=0$ or $c_j' = 0$,}\end{cases}\end{equation}
where $\alpha$ is the solution in $\lbrack 0,1\rbrack$ to
\begin{equation}\label{eq:otrol}
\mathop{\sum_{j\in J}}_{c_j,c_j'\ne 0} c_j^{\alpha} (c_j')^{1-\alpha} 
\left(\frac{d-1}{d} - \left(\log \frac{c_j}{c_j'}\right) 
\left(\frac{1}{d}+\frac{d-1}{d} \alpha \right) \right) = \theta ,\end{equation}
if there is a solution in $\lbrack 0,1\rbrack$ (in which case
it is unique). If (\ref{eq:otrol}) has no solution in $\lbrack 0,1\rbrack$,
then $\alpha$ is either $0$ or $1$, depending on which of the two resulting
choices of $\vec{x}$ (as per (\ref{eq:horna})) gives the smaller value of
$g(\vec{x})$.
If the sum in (\ref{eq:otrol}) has no terms, then $\vec{x}$ is the zero vector.
\end{lem}
When (\ref{eq:otrol}) has no solutions in $\lbrack 0,1\rbrack$, 
the minimal value of $g(\vec{x})$ is
easy to describe: as we shall see, it equals 
\begin{equation}\label{eq:luverly}
\max\left(1 - \frac{1}{d} \sum_{j\in J_0} c_j, \theta + 1 - \sum_{j\in J_0}
c_j'\right),\end{equation}
 where $J_0 = \{j\in J : c_j, c_j' \ne 0\}$. 

No matter
whether (\ref{eq:otrol})
has a solution in $\lbrack 0,1\rbrack$ or not, the minimal value of
$g(\vec{x})$ will be greater than $\theta$ if and only if 
$\frac{d-1}{d} + \frac{1}{d} I_{\vec{c}}(\vec{c}') > \theta$. This is 
easy to see: if $\frac{d-1}{d} + \frac{1}{d} I_{\vec{c}}(\vec{c}')\leq \theta$,
then $g(\vec{c}') \leq \theta$, and so $\min g(\vec{x}) \leq \theta$;
if $\frac{d-1}{d} + \frac{1}{d} I_{\vec{c}}(\vec{c}') > \theta$, we have
$\frac{d-1}{d} + \frac{1}{d} I_{\vec{c}}(\vec{x}) > \theta$ for
$\vec{x}$ in a neighbourhood of $\vec{c}'$, and $\theta + I_{\vec{c}'}(\vec{x})
> \theta$ outside the neighbourhood.
\begin{proof}[Proof of Lemma \ref{lem:mamsol}]
For every $j\in J$ such that $c_j=0$ or $c_j'=0$, the variable $x_j$
is forced to be zero for all $\vec{x}\in X$ such that 
$g(\vec{x}) < \infty$. At the same time, if $x_j=0$, the terms involving
$x_j$ make no contribution\footnote{We set the 
convention $0\log 0 = 0$ when $I_{\vec{c}}(\vec{x})$ was defined.
See the comment after (\ref{eq:jacbo}).}
 to
either $I_{\vec{c}}(\vec{x}) = 1 - \sum_j x_j + \sum_j x_j \log 
\frac{x_j}{c_j}$ or $I_{\vec{c}'}(\vec{x}) = 1 - \sum_j x_j + \sum_j x_j \log 
\frac{x_j}{c_j'}$. (See the convention on $0 \log 0$ chosen after
(\ref{eq:jacbo}).) Hence, we may redefine $J$ to be $\{j\in J:c_j,c_j'\ne 0
\}$, and thus reduce the problem to the case in which $c_j, c_j' \ne 0$
for all $j\in J$. We assume, then, that $c_j, c_j'\ne 0$ for all $j\in J$;
consequently, $I_{\vec{c}}(\vec{x})$ and $I_{\vec{c}'}(\vec{x})$ will be
smooth on $(\mathbb{R}^+)^J$, which is an open superset of $X$.

Define $g_1(\vec{x}) = \frac{d-1}{d} + \frac{1}{d} I_{\vec{c}}(\vec{x})$,
$g_2(\vec{x}) = \theta + I_{\vec{c}'}(\vec{x})$.  Then 
$g(x)$ equals $\max(g_1(x),g_2(x))$.
If $g(\vec{x})$ is minimal on $X$ at $\vec{x}\in X$, then it is minimal on
$(\mathbb{R}^+)^J$ at $\vec{x}$: the partial derivatives
$\frac{\partial}{\partial x_j} g_1(\vec{x})$, $\frac{\partial}{\partial x_j}
g_2(\vec{x})$ are negative for $x_j < \min(c_j,c_j')$ and positive for
$x_j > \max(c_j,c_j')$. 
 The minimum of $g$ on $(\mathbb{R}^+)^J$
 may be
attained at a point $\vec{x}$ where
\begin{enumerate}
\item\label{it:law1} $g_1(\vec{x})$ has a local minimum,
\item\label{it:law2} $g_2(\vec{x})$ has a local minimum, or
\item\label{it:dere} $g_1(\vec{x}) = g_2(\vec{x})$.
\end{enumerate}
(There are no other cases: if none of the above were to
hold, a small displacement in
$\vec{x}$ will decrease whichever one of $g_1(\vec{x})$ or $g_2(\vec{x})$
is greater, and thereby decrease $g(\vec{x})=\max(g_1(\vec{x}),g_2(\vec{x}))$
from its supposed minimum.) The only local
 minimum of $g_1(\vec{x})$ on $(\mathbb{R}^+)^J$ is at
$\vec{x} = \vec{c}$, and the only local maximum of $g_2(\vec{x})$ on 
$(\mathbb{R}^+)^J$ is at
$\vec{x} = \vec{c}'$. It remains to consider case (\ref{it:dere}).
Then $g$ reaches a minimum on $(\mathbb{R}^+)^J$ at a point $\vec{x}$
on the surface
$S$ described by the equation 
$g_1(\vec{x}) = g_2(\vec{x})$. By restriction, $g$ reaches a minimum
on $S \cap (\mathbb{R}^+)^J$ at $\vec{x}$.
Now, on $S\cap (\mathbb{R}^+)^J$, the function $g_1(\vec{x})$ equals
$g(\vec{x})$. It follows that 
$\nabla g_1(\vec{x})$ is perpendicular to $S$, and thus $\nabla g_1(\vec{x})$
is a scalar 
multiple of $\nabla (g_1(\vec{x}) - g_2(\vec{x}))$. In other words,
one of $\nabla g_1(\vec{x})$, $\nabla g_2(\vec{x})$ is a scalar multiple of
the other. Now
\[\begin{aligned}
\nabla g_1(\vec{x}) &= \nabla \left(\frac{d-1}{d} + \frac{1}{d} 
I_{\vec{c}}(\vec{x})\right) = \frac{1}{d} \nabla I_{\vec{c}}(\vec{x}) =
\left\{\frac{1}{d} \log \frac{x_j}{c_j} \right\}_{j\in J}\\
\nabla g_2(\vec{x}) &= \nabla (\theta + I_{\vec{c}'}(\vec{x})) = 
\nabla I_{\vec{c}'}(\vec{x}) = \left\{\log \frac{x_j}{c_j'}\right\}_{j\in J} .
\end{aligned}\]
We conclude that $x_j = c_j^{\alpha} (c_j')^{1-\alpha}$ for some $\alpha$.
Since $\vec{x}\in X$, we know that $\alpha$ must be in $\lbrack 0,1\rbrack$.
As we have already seen,
$x_j = c_j^{\alpha} (c_j')^{1-\alpha}$ holds
in cases (\ref{it:law1}) and (\ref{it:law2}) just
 as well, with $\alpha = 1$ and $\alpha=0$, respectively.

Our task is now to find the minimum of 
\[g(\{c_j^{\alpha} (c_j')^{1-\alpha}\}) = \max(g_1(
\{c_j^{\alpha} (c_j')^{1-\alpha}\}),g_2(
\{c_j^{\alpha} (c_j')^{1-\alpha}\}))\]
for $\alpha\in \lbrack 0,1\rbrack$. The map $h_1:\alpha \mapsto
g_1(\{c_j^{\alpha} (c_j')^{1-\alpha}\})$ is increasing, whereas
$h_2:\alpha \mapsto g_2(\{c_j^{\alpha} (c_j')^{1-\alpha}\})$ 
is decreasing. Thus,
\begin{equation}\label{eq:otror}
g_1(\{c_j^{\alpha} (c_j')^{1-\alpha}\}) =
 g_2(\{c_j^{\alpha} (c_j')^{1-\alpha}\})\end{equation}
 for at most one $\alpha\in \lbrack
0,1\rbrack$, and, if such an $\alpha$ exists, $g(\{c_j^{\alpha} 
(c_j')^{1-\alpha}\})$ attains its minimum thereat. 
Writing out $g_1$ and $g_2$, we see that (\ref{eq:otror}) is equivalent
to (\ref{eq:otrol}). 

If (\ref{eq:otror}) has no solution $\alpha$ within $\lbrack 0,1\rbrack$,
the minimum of
$g(\{c_j^{\alpha} 
(c_j')^{1-\alpha}\})$ is $\min(\max(h_1(0),h_2(0)),
\max(h_1(1),h_2(1)))$, which, since $h_1$ is increasing
and $h_2$ is decreasing, equals $\max(h_1(0),h_2(1))$, which, written
in full,
is
\[
\max\left(\frac{d-1}{d} + \frac{1}{d} \left(1 - \sum_{j\in J} c_j\right),\;
\theta + 1 - \sum_{j\in J} c_j'\right) .\]
This is nothing other than (\ref{eq:luverly}).
\end{proof}

Thanks to Prop.\ \ref{prop:knut}, Lem.\ \ref{lem:knotty} 
and Lem.\ \ref{lem:mamsol}, we finally know how to bound the number
of elements $n\in S$ ($S$ a tight set) such that $p^k|f(n)$ for
some large prime $p$. Our end is to estimate the number of elements 
$n\in S$ such that $f(n)$ is free of $k$th powers. The remaining way
to the end is rather short.
\begin{lem}\label{lem:sourc}
Let $f\in \mathbb{Z}\lbrack x\rbrack$ be a polynomial. Let $k\geq 2$. Let
$S$ be a predictable, $(P,\theta)$-tight set. Then the number of elements 
$n\in S\cap \lbrack 1,N\rbrack$ such that $f(n)$ is free
of $k$th powers equals
\begin{equation}\label{eq:sarko}\begin{aligned}
\sum_{d\leq D} \mu(d) &|\{n\in S\cap \lbrack 1, N\rbrack : d^k |f(n)\}|
+ O_f(D^2 + D^{-(k-1)} (\log D)^c N)\\
&+ O(|\{n\in S\cap \lbrack 1,N\rbrack: \exists p>D^2 \text{ \st\; }
p^k|f(n)\}|)\end{aligned}\end{equation}
for every $D\geq 2$ and some $c>0$ depending only on $\deg(f)$. The
second implied constant is absolute.
\end{lem}
Cf.\ \cite{Hsq}, Prop.\ 3.4.
\begin{proof}
Apply the {\em riddle} in \cite{Hsq}, \S 3 (that is, \cite{Hsq}, Prop.\
3.2) with $\mathscr{P}$ equal to the set of all primes,
$\mathscr{A} = S\cap \lbrack 1,N\rbrack$, $r(a) = \{\text{$p$ prime}:
p^k|a\}$, $f(a,d)=1$ if $d=\emptyset$, and $f(a,d)=0$ for $d$ non-empty.
Use the bound 
\[|\{n\in S\cap \lbrack 1,N\rbrack : d^k | f(n)\}| \leq
\{1\leq n\leq N : d^k|f(n)\} \ll_f
\frac{N (\deg f)^{\omega(d)}}{d^k} + 1.\]
\end{proof}

Now it only remains to reap the fruits of our labour.
In order to avoid unnecessarily lengthy and complicated statements, we will 
give the explicit results below for irreducible polynomials $f$ alone.
They can be easily restated for general polynomials in the manner of
the statement of the main theorem.
\begin{prop}\label{prop:ints}
Let $f\in \mathbb{Z}\lbrack x\rbrack$ be a polynomial. Let 
$k\geq \max(2,\deg(f)-1)$. 
If $k = \deg(f) - 1$, then
\begin{equation}\label{eq:antar}\begin{aligned}
|\{1\leq n\leq N: d^k|f(n) \Rightarrow d=1\}| 
 &=  
N\cdot \prod_p \left(1 - \frac{\rho_f(p^k)}{p^k}\right) \\ &+ 
O_{f,k,\epsilon}\left(N (\log N)^{- \left(1 - 
\frac{\sigma(\Gal_f)}{\deg(f) |\Gal_f|}\right)
 + \epsilon}\right)\end{aligned}\end{equation}
for every $\epsilon>0$, where $\rho_f(p^k)$ is the number of solutions
to $f(x)\equiv 0 \mo p^k$ in $\mathbb{Z}/p^k \mathbb{Z}$
and $\sigma(\Gal_f)$ is the number of maps in $\Gal_f$ that 
have fixed points.

If $k\geq \deg(f)$, then
\begin{equation}\label{eq:ogoh}\begin{aligned}
|\{1\leq n\leq N: d^k|f(n) \Rightarrow d=1\}| 
&= 
N\cdot \prod_p \left(1 - \frac{\rho_f(p^k)}{p^k}\right) \\ &+ 
O_{f,k}\left(N^{1 - \delta}\right)\end{aligned}\end{equation}
for some $\delta$ depending only on $\deg(f)$.
\end{prop}
\begin{remark}
The value of $\delta$
in the right side of 
(\ref{eq:ogoh}) may be chosen to be
$1- \frac{\max(2,\deg(f))}{k+1} - \epsilon'$, with
$\epsilon'>0$ arbitrarily small: sharpen
Lem.\ \ref{lem:knotty} as detailed in the remark after its proof,
and then apply Lemmas \ref{lem:knotty} and \ref{lem:sourc}
with $D = N^{1/(k+1)}$. 
\end{remark}
\begin{proof}
Apply Lemma \ref{lem:sourc} with $D = N^{\delta'}$, where $\delta'>0$
will be chosen later. The error term 
$|\{n\in S\cap \lbrack 1, N\rbrack : \exists p>D^2 \text{ \st\; } p^k|f(n)\}|$
in (\ref{eq:sarko}) is at most $N$ times (\ref{eq:schuby}), and thus
can be bounded by Prop.\ \ref{prop:knut}. The exponent $\gamma$ in
the bound on (\ref{eq:schuby}) in Prop.\ \ref{prop:knut} can be determined
by Lemma \ref{lem:mamsol}; it amounts to $1 - \frac{\sigma(\Gal_f)}{
\deg(f) |\Gal_f|}$. We are left with the main term
$\sum_{d\leq N^{\delta'}} \mu(d)\cdot  |\{1\leq n\leq N :
d^k |f(n)\}|$; we wish to show that it equals $N\cdot
\prod_p \left(1 - \frac{\rho_f(p^k)}{p^k}\right)$ plus a small error term. 

Since
$|\{1\leq n\leq N :
d^k |n\}| = N/d^k + O(1)$, we have
\[\sum_{d\leq N^{\delta'}} \mu(d)\cdot |\{1\leq n\leq N :
d^k |f(n)\}| = \sum_{d\leq N^{\delta'}} \mu(d) \rho_f(d^k)
\left(\frac{N}{d^k} + O(1)\right) ,\]
where $\rho_f(m) = \prod_{p|m} \rho_f\left(p^{v_p(m)}\right)$.
The right side equals
\begin{equation}\label{eq:indque}
N \sum_{d} \mu(d) \frac{\rho_f(d^k)}{d^k} + 
O\left(\sum_{d\leq N^{\delta'}} |\mu(d) \rho_f(d^k)| +
N \sum_{d>N^{\delta'}} \frac{|\mu(d) \rho_f(d^k)|}{d^k} \right).\end{equation}
By $|\mu(d) \rho_f(d^k)| \leq \prod_{p|d} |\rho_f(p^k)|\ll_f
(\deg f)^{\omega(d)}$, we have that (\ref{eq:indque}) equals
$N\cdot \prod_p \left(1 - \frac{\rho_f(p^k)}{p^k}\right)$ plus \[ 
O_f\left(N^{\delta'} (\log N)^{\deg(f)-1}
+ N^{1 - (k-1)\delta'} (\log N)^{\deg(f)-1}\right) \]
plus the error term $O(N^{2 \delta'} + N^{-(k-1) \delta'} (\log D)^c N)$
coming from Lemma \ref{lem:sourc}.
We set $\delta' = \frac{1}{k}$ and are done.
\end{proof}
\begin{prop}\label{prop:primos}
Let $f\in \mathbb{Z}\lbrack x\rbrack$ be a polynomial. Let 
$k\geq \max(2,\deg(f)-1)$. If $k = \deg(f) - 1$, then
\begin{equation}\label{eq:oglor}\begin{aligned}
|\{\text{$q$ prime, $q\leq N$}: d^k|f(q) \Rightarrow d=1
\}| &= 
\pi(N)\cdot \prod_p \left(1 - \frac{\rho_{f,*}(p^k)}{p^k - p^{k-1}}\right) 
\\ &+ 
\pi(N)\cdot O_{f,k,\epsilon}\left((\log N)^{-\gamma + \epsilon}\right)
\end{aligned}\end{equation}
for every $\epsilon>0$, where $\pi(N)$ is the number of primes from
$1$ to $N$, $\rho_{f,*}(p^k)$ is the number of solutions
to $f(x)\equiv 0 \mo p^k$ in $(\mathbb{Z}/p^k \mathbb{Z})^*$,
and $\gamma = g(\vec{x})-1$, where $g$ is as in (\ref{eq:empan}) and
$\vec{x}$ is as in (\ref{eq:horna}), with $\vec{c}$ and $\vec{c}'$
as in Prop.\ \ref{prop:knut}. We have $\gamma>0$ if and only if
$I_f>1$.

If $k\geq \deg(f)$, then, for every $A>0$,
\begin{equation}\label{eq:anda}\begin{aligned}
|\{\text{$q$ prime, $q\leq N$}: d^k|f(q) \Rightarrow d=1
\}| &= 
\pi(N)\cdot \prod_p \left(1 - \frac{\rho_{f,*}(p^k)}{p^k - p^{k-1}}\right) \\ &+ 
O_{f,k,A}\left(N (\log N)^{-A}\right). \end{aligned}\end{equation}
\end{prop}
\begin{proof}
Proceed as in the proof of Prop.\ \ref{prop:ints}, with 
$D = (\log N)^{A}$; use Siegel-Walfisz to estimate
$|\{\text{$q$ prime}, q\leq N: d^k|f(q)\}|$ for $d\leq D$.
\end{proof}
\begin{remark}
 Weaker effective results can
be used instead of Siegel-Walfisz;
if the best such results are used (see, e.g., \cite{Da}, \S 14, (9),
and \S 20, (11)) then, as can be shown by a simple computation,
the error term in (\ref{eq:oglor}) remains unaltered 
for $\deg(f)\leq 6$.
\end{remark}
\begin{remark}
See Table \ref{tab:bogo} for the values of $\gamma$ for $\deg(f)\leq 6$.
If $\deg(f)>6$, Nair's result (\cite{Na}, Thm.\ 3) applies; its error
term is no larger than $O_A(N (\log N)^{-A})$, where $A>0$ is arbitrarily
large.
\end{remark}
\begin{remark}
Let $f$ be irreducible and {\em normal}; that is,
assume its degree $d$ equals the degree $\left|\Gal_f\right|$ of its
splitting field. 
As seen in (\ref{eq:squiggle}),
we have $I_f>1$, and so $\gamma>0$ in 
(\ref{eq:oglor}); in other words, the error term 
is smaller than the main term. Because
the structure of $\Gal_f$ is particularly simple, we shall be able to
give a fairly uncomplicated expression for $\gamma$.

For $x\in (-e^{-1},0)$, let $W_{-1}(x)$ be the smaller of the two solutions
$y$ to $y e^y = x$. (As can be gathered from the notation, $W_{-1}$ is
one of the branches $W_k$ of the Lambert $W$ function.) Let $\vec{c}$,
$\vec{c}'$ be as in Prop.\ \ref{prop:knut}. By the Cauchy-Frobenius formula
and the fact that $f$ is normal, the only map in $\Gal_f$ with any fixed
points is the identity. Thus $c_{\langle g \rangle}' = 0$
for $g\ne e$ and $c_{\langle e\rangle}' = 1$, while $c_{\langle e \rangle}
= \frac{1}{d}$.  Therefore, (\ref{eq:otrol}) can be rewritten as
\begin{equation}\label{eq:munic}
\left(\frac{1}{d}\right)^{\alpha} \left(\frac{d-1}{d} + \log(d) \left(
\frac{1}{d} + \frac{d-1}{d} \alpha\right)\right) = 1 .\end{equation}
(We have $\theta=1$ because we are working on the primes.)
We let $y = -\log(d) \alpha - \left(1 + \frac{\log d}{d-1}\right)$ and rewrite
(\ref{eq:munic}) as
\begin{equation}\label{eq:bana}
y e^y = \frac{- d^{(d-2)/(d-1)}}{e (d-1)} .\end{equation}
Since $d-1 > d^{(d-2)/(d-1)}$ for $d\geq 3$ (as is our case),
the right side of (\ref{eq:bana}) is in the range
$(-e^{-1},0)$, and thus
 $y = W_{-1}\left(\frac{- d^{(d-2)/(d-1)}}{e (d-1)}\right)$.
Hence
\[\alpha = -\frac{1}{\log(d)} \left( 
W_{-1}\left(\frac{- d^{(d-2)/(d-1)}}{e (d-1)}\right) + 1 + \frac{\log d}{
d-1}\right) .\] Now (\ref{eq:horna}) gives $x_{\langle e\rangle} = 
d^{-\alpha}$, $x_{\langle g\rangle} = 0$ for $g\ne e$
and (\ref{eq:empan}) yields \[\begin{aligned}
\gamma &= g(\vec{x}) - 1 = I_{\vec{c}'}(\vec{x}) = 
1 - \frac{1 + \alpha \log d}{d^{\alpha}} \\ &=
- \frac{d \log d}{(d-1)^2 W_{-1}\left(\frac{- d^{(d-2)/(d-1)}}{e (d-1)}\right)}
- \frac{1}{d-1}.\end{aligned}\]
Thus,
for $d$ large, $\gamma \sim \frac{d \log d}{(d-1)^2} - \frac{1}{d-1}$,
which goes to $0$ as $d\to \infty$.
\end{remark}
\begin{cor}[to Prop.\ \ref{prop:primos}]
Let $f\in \mathbb{Z}\lbrack x\rbrack$ be a cubic polynomial irreducible
over $\mathbb{Q}$. Suppose that its discriminant is a square. Then the number
of primes $q$ from $1$ to $N$ such that $f(q)$ is square-free equals
\[\pi(N) \cdot \prod_p \left(1 - \frac{\rho_{f,*}(p^2)}{p^2 - p}\right) +
O_\epsilon(\pi(n) \cdot (\log N)^{-\gamma + \epsilon})\]
for every $\epsilon>0$, 
where 
\begin{enumerate}
\item $\pi(N) = N/\log N + O(N/(\log N)^2)$ is the number of primes 
up to $N$,
\item $\rho_{f,*}(p^2)$ is the number of solutions to $f(x)\equiv 0 \mo p^2$
in $(\mathbb{Z}/ p^2 \mathbb{Z})^*$, 
\item $\gamma$ equals $1 - 3^{-\alpha} + 3^{-\alpha} \log 3^{-\alpha}>0$,
where $\alpha$ is the only solution in $\lbrack 0,1\rbrack$ to 
$3^{-\alpha} \left(\frac{2}{3} - \left(\log \frac{1}{3}\right) \cdot
\left(\frac{1}{3} + \frac{2}{3} \alpha\right)\right) = 1$.
\end{enumerate}
\end{cor}
Numerically, $\gamma = 0.003567\dotsc$.
\begin{proof}
Since the discriminant of $f$ is a square, the Galois group $\Gal_f$
of $f$ is $A_3$. Apply Prop.\ \ref{prop:primos}.
\end{proof}
\begin{table}
\begin{center}
\begin{tabular}{ll|ll|ll}
$\Gal_f$ & $\gamma$ 
& $\Gal_f$ & $\gamma$ 
& $\Gal_f$ & $\gamma$\\ \hline
$A_3$ &   $0.0035671$ 
& & 
& &
\\\hline 
$C(4)$ & $0.0265166$
& $E(4)$ & $0.0265166$
& $D(4)$ & $0.0006060$
\\\hline
$C(5)$ & $0.0417891$ & & & & \\\hline
$C(6)$ & $0.0505865$
& $D_6(6)$ & $0.0505865$ 
& $D(6)$ & $0.0104233$\\
$A_4(6)$ & $0.0104233$ & 
$F_{18}(6)$ & $0.0170657$ & $2A_4(6)$ & $0.0157592$ \\
$F_{18}(6):2$ & $0.0000529$ & 
$F_{36}(6)$ & $0.0000529$ & $2S_4(6)$ & $0.0000059$\\
\\
\end{tabular}
\caption{Values of
$\gamma$
for $\theta=1$ and $f$ of degree $d\leq 6$ with $I_f>1$. The number of primes $p\leq N$ with $f(p)$ free
of $(d-1)$th powers equals a constant $c_f$ times $\pi(N)
(1 + O((\log N)^{-\gamma+ \epsilon}))$.}\label{tab:bogo}
\end{center}
\end{table}
\begin{table}
\begin{center}
\begin{tabular}{ll|ll|ll}
$\Gal_f$ & $\gamma$ 
& $\Gal_f$ & $\gamma$ 
& $\Gal_f$ & $\gamma$\\ \hline
$A_3$ &   $0.3888889$ 
& $S_3$ & $0.2777778$
& &
\\\hline 
$C(4)$ & $0.4375000$
& $E(4)$ & $0.4375000$ & $D(4)$ & $0.4062500$\\
$A_4$ & $0.3125000$ & $S_4$ & $0.3437500$ & &\\\hline
$C(5)$ & $0.4600639$
& $D(5)$ & $0.3800000$ & $F(5)$ & $0.3400000$\\
$A_5$ & $0.3800000$ & $S_5$ & $0.3733333$ & &\\\hline
$C(6) $ & $0.4728484$
& $D_6(6) $ & $0.4728484$ & $D(6) $ & $0.4444444$\\
$A_4(6) $ & $0.4444444$ & $F_{18}(6) $ & $0.4537037$
&$2A_4(6) $ & $0.4513889$ \\$S_4(6d) $ & $0.4305556$
&$S_4(6c) $ & $0.4305556$ & $F_{18}(6):2 $ & $0.4351852$\\
$F_{36}(6) $ & $0.4351852$ & $2S_4(6) $ & $0.4340278$
&$L(6) $ & $0.3888889$\\
$F_{36}(6):2 $ & $0.4259259$ &$L(6):2 $ & $0.4027778$ & 
$A_6 $ & $0.3935185$\\
$S_6 $ & $0.3946759$ & & & &
\\
\end{tabular}
\caption{
Values of $\gamma$
for $\theta=\frac{1}{2}$ and $f$ of degree $d\leq 6$. The number of 
sums of two squares $q\leq N$ with $f(q)$ free
of $(d-1)$th powers equals a constant $c_f$ times $\varpi(N) (1 + O((\log N)^{-\gamma+ \epsilon}))$, where $\varpi(N)$ is the sum of integers up to $N$
that can be written as sums of two squares.}\label{tab:mipi}
\end{center}
\end{table}
\begin{prop}
Let $f\in \mathbb{Z}\lbrack x\rbrack$ be a polynomial. Let 
$k\geq \max(2,\deg(f)-1)$. Let $S$ be the set of all integers that
are the sum of two squares.
If $k = \deg(f) - 1$, then
\begin{equation}\label{eq:osto}\begin{aligned}
|\{n\in S\cap \lbrack 1,N\rbrack: d^k|f(n) \Rightarrow d=1
\}| &= 
 \varpi(N) \cdot \prod_p (1 - \rho_{f,\circ}(p^k)) \\ &+
O_{f,k,\epsilon}\left(\varpi(N) \cdot (\log N)^{-\gamma + \epsilon}\right),
\end{aligned}\end{equation}
for every $\epsilon>0$, where 
\[\begin{aligned}
\varpi(N) &= |S\cap \lbrack 1, N\rbrack| \sim 
\left(2\cdot  \prod_{p\equiv 3 \mo 4} (1 - p^{-2})\right)^{-1/2} \cdot
\frac{N}{\sqrt{\log N}}
,\\
\rho_{f,\circ}(p^k) &= 
\mathop{\sum_{a\in \mathbb{Z}/p^k \mathbb{Z}}}_{f(a)\equiv 0 \mo p^k} 
\rho_{\circ}(a,p^k),\\
\rho_{\circ}(a,p^k) &= 
\begin{cases}
p^{-k} (1 + p^{-1}) &\text{if $p\equiv 3 \mo 4$, $v_p(a)$ even,
$v_p(a)<k$,}\\
0 &\text{if $p\equiv 3 \mo 4$, $v_p(a)$ odd, $v_p(a)<k$,}\\
p^{-k} &\text{if $p\equiv 3 \mo 4$, $v_p(a)$ even, $v_p(a)=k$,}\\
p^{-(k+1)} &\text{if $p\equiv 3 \mo 4$, $v_p(a)$ odd, $v_p(a)=k$,}\\
p^{-k} &\text{otherwise,}
\end{cases}\end{aligned}\]
and $\gamma = g(\vec{x}) - \frac{1}{2}$, 
where $g$ is as in (\ref{eq:empan}) and
$\vec{x}$ is as in (\ref{eq:horna}), with $\vec{c}$ and $\vec{c}'$
as in Prop.\ \ref{prop:knut}.

If $k \geq \deg(f)$, then, for all $A>0$,
\begin{equation}\begin{aligned}
|\{n\in S\cap \lbrack 1,N\rbrack: d^k|f(n) \Rightarrow d=1
\}| &= \varpi(N) \cdot \prod_p (1 - \rho_{f,\circ}(p^k)) \\ &+
O_{f,k,A}\left(N (\log N)^{-A}\right),
\end{aligned}\end{equation}
where $\varpi(N)$ and $\rho_{f,\circ}$ are as above.
\end{prop}
\begin{proof}
Proceed as in the proof of Prop.\ \ref{prop:ints}. Use \cite{Ri},
Hilfs\"atze 10 and 12, to show 
\begin{equation}\label{eq:otrumo}
|S\cap \lbrack 1,N\rbrack \cap (a + m \mathbb{Z})| =
\rho_{\circ}(a,m) \cdot |S \cap \lbrack 1,N\rbrack| + O_A(N e^{-c \sqrt{\log N}})
\end{equation}
for $a$, $m$ with $\gcd(a,2 m)=1$, where $\rho_{\circ}(a,m) =
\prod_{p|m} \rho_{\circ}(a,p^{v_p(m)})$ and $c$ is a positive constant.
Extend (\ref{eq:otrumo}) to the case $\gcd(a,2 m) \ne 1$ by direct use of the
fact that $n\in S$ if and only if $v_p(n)$ is even for every $p\equiv 3 \mo 4$.
\end{proof}
\begin{remark}
Since $\frac{d-1}{d} + \frac{1}{d} I_{\vec{c}}(\vec{x}) \geq
\frac{d-1}{d} \geq \frac{2}{3}$ for all $\vec{x}$, we have
$\gamma>\theta$ for $\theta = 0.5$ and $f$ arbitrary. Hence,
the error term in (\ref{eq:osto}) is smaller than
$O(N (\log N)^{-1/2 + \epsilon})$, $\epsilon>0$ arbitrary, and thus it is
smaller than the main term
for all $f$ for which the infinite product in (\ref{eq:osto})
does not vanish. 

The values of $\gamma$ for $\deg(f)\leq 6$ are in Table \ref{tab:mipi}.
Most entries in the table are rational; this is because,
when $\theta = \frac{1}{2}$ and $\vec{c}$, $\vec{c}'$ are as in
Prop.\ \ref{prop:knut}, the minimum of (\ref{eq:empan}) is
reached at $\vec{x} = \vec{c}'$ for many (but not all) $f$.
\end{remark}
\begin{center}
* * *
\end{center}

The approach taken in this paper can be applied to improve upon
the error term given in \cite{Hsq} for the estimated number of
pairs of integers $(x,y)\in \lbrack 1, N\rbrack^2$ such that
$f(x,y)$ is square-free, where $f$ is a sextic homogeneous polynomial.
Asymptotics were first given in \cite{Gre}, with the error
term $O(N^2 (\log N)^{-1/3})$; soon thereafter, K. Ramsay \cite{Ra}
attained $O(N^2 (\log N)^{-1/2})$ by means of a slight modification
in the argument. The error term in \cite{Hsq} depends on the
Galois group of $f(x,1)$; for $f$ generic, it is 
$O(N^2 (\log N)^{-0.7043\dotsc})$. We can do better now for every $f$.
\begin{prop}\label{prop:azumo}
Let $f\in \mathbb{Z}\lbrack x,y\rbrack$ be a homogeneous
polynomial of degree $6$ irreducible in 
$\mathbb{Q}\lbrack x,y\rbrack$.
Then the number of pairs of
integers $(x,y)$, $1\leq x,y \leq N$, such that 
$f(x,y)$ is square-free equals
\begin{equation}
N^2 \prod_{p} \left(1 - \frac{\rho_f(p^2)}{p^4}\right) +
O_{f,\epsilon}\left(N^2 (\log N)^{-1 + \frac{\sigma(G)}{3 |G|} + \epsilon}\right),
\end{equation}
where $G$ is the Galois group of the splitting field of $f_0(x)
= f(x,1)$, 
$\sigma(G)$ is the number of maps $g$ in $G$ that have fixed
points,
$\rho_f(p^2)$ is the number of solutions 
$(x,y)\in (\mathbb{Z}/p^2 \mathbb{Z}) \times
(\mathbb{Z}/p^2 \mathbb{Z})$
to 
 $f(x,y)\equiv 0 \mo p^2$, and $\epsilon>0$ is arbitrary. The implied
constant depends only on $f$ and $\epsilon$.
\end{prop}
\begin{proof}[Proof (Sketch)]
Proceed as in Prop.\ \ref{prop:ints}, replacing
Lem.\ \ref{lem:sourc} by \cite{Hsq},
Prop.\ 3.5. It remains to bound
\begin{equation}\label{eq:prior}
|\{1\leq x,y\leq N: \exists p>N^2 \text{\; s.t. } p^2 |f(x,y)\}| .
\end{equation}
This we do by giving a bound for the number of rational points
on $C_r : r y^2 = g(x)$ with $1\leq r\leq N^2$, where $g(x) = f(x,1)$. 
We can do this by finding for the great
majority of $r$ (as in Prop.\ \ref{prop:notung})
a divisor $t_0|r$, $t_0 > N^{2-\epsilon}$, with few prime factors,
and then using it as 
in \cite{HV}, \S 5.
(As before, a divisor $t_0|r$ of the right size will exist for all $r$
outside a small set, viz., a set of cardinality $\ll N^2 (\log N)^{-A}$,
$A$ arbitrary.)
The divisor $t_0$ is large enough to increase the angle given by
Mumford's gap principle to $\pi/2 - O(\epsilon)$. We can then
apply sphere-packing results (Lem.\ \ref{lem:sphpack}), bounding the
rank of $C_r$ in terms of $w(r)$ as in \cite{Hsq}, Prop.\ 4.22
(that is, using \cite{Ca}, though we may use the more general
statements in \cite{PSc} instead). Let $D$ be a positive integer that
will be set later.
We consider all $r\leq D$
such that (a) $r$ has a divisor $t_0$ as above, and
(b) the 
Frobenius element in $G = \Gal_{g}$ of every $p|r$ has fixed points.
There are 
$\ll_f D (\log N)^{- 1 + \frac{\sigma(G)}{|G|}}$
such integers $r$. Since the bound on the number of rational
points per $r$ coming from sphere-packing is $(\log N)^\epsilon$
for all $r$ outside a small set, we obtain a total bound of
$\ll_f D (\log N)^{- 1 + \frac{\sigma(G)}{|G|} + \epsilon}$.
Since $r p^2 = f(x)$ with $r>D$ implies $p\ll N^3/\sqrt{D}$, we can
bound the contribution to (\ref{eq:prior}) of solutions to
$r p^2 = f(x)$ with $r>D$ by $O_f\left(\frac{N^3}{(\log N) \sqrt{D}}\right)$.
Set $D = (\log N)^{-\frac{2 \sigma(G)}{3 |G|}} N^2$. We conclude that
(\ref{eq:prior}) is at most $O_f(N^2 (\log N)^{-1 + \frac{\sigma(G)}{3 |G|}
+\epsilon})$.
\end{proof}
\begin{remark}
Pairs
of integers $(x,y)\in \mathbb{Z}^2 \cap \lbrack 1,N\rbrack$ are numerous
enough that considerations of entropy are not needed to prove
Prop.\ \ref{prop:azumo}, and, in fact, would not help.
Thus, the situation is similar to that in Prop.\ \ref{prop:ints},
and the contrary of the situation in every other result in this paper.
\end{remark}

\end{document}